# First-order natural deduction in Agda

Louis Warren
School of Mathematics and Statistics, University of Canterbury, New Zealand

## Contents

# 1 Introduction

Agda [2] is a dependently-typed functional programming language [9], based on an extension of intuitionistic Martin-Löf type theory. We implement first order natural deduction in Agda. We use Agda's type checker to verify the correctness natural deduction proofs, and also prove properties of natural deduction, using Agda's proof assistant functionality. This implementation corresponds to a formalisation of natural deduction in constructive type theory, and the proofs are verified by Agda to be correct (under the assumption that Agda itself is correct).

The Agda code below has been written in *literate Agda*, which allows Agda to be mixed with LaTeX. The files which have been used to typeset this document can also be evaluated and type checked. Some trivial proofs are omitted from the typeset document; these are only hidden for brevity, and are still present in the code and used by Agda. The results which rely on them are therefore still verified. This should not be mistaken for use of postulates, wherein Agda itself is told to assume that a proof exists. Postulates are used only in the module for outputting natural deduction proofs as LaTeXfor use with the *bussproofs* package. All other code type checks with Agda in safe mode, meaning that it provably halts.

Each of the following sections corresponds to its own literate Agda file. Sections named with a file name ending in '.lagda' are modules. Each section imports the modules preceding it, unless stated otherwise. These module declarations and imports have been hidden for brevity.

Inspiration for the definition of vector types and decidable types comes from the Agda standard library [1]. However, the standard library will not be directly imported, to maintain clarity of definitions, and because it is not needed. We will use built-in types for natural numbers, lists, and the dependent sum, explaining their definitions when they appear.

The full code is available online at `https://lsw.nz/tome`. It is written for (and has been tested against) Agda version 2.6.0.

# 2 Decidable.lagda

We begin with a module which defines decidability.

Agda has a built-in module defining equality. We import this module and re-export it here. For illustrative purposes, a simplified version of this definition for small types (types of type `Set`) is commented below.

```
open import Agda.Builtin.Equality public

{-
  data _≡_ {A : Set} (x : A) : A → Set where
    refl : x ≡ x
-}
```

For every $x$ of any type, there is a constructor for $x \equiv x$. An instance of the equality $x \equiv y$ is a proof that $x$ and $y$ are intensionally equal. In Agda, we use `data` types as a convenient notation for what would otherwise be defined type-theoretically using W-types.

The bottom type, ⊥, has no constructors, and so is provable only from absurdity. The usual definition of negation follows, as does an abbreviation for inequality.

```
data ⊥ : Set where

¬_ : (A : Set) → Set
¬ A = A → ⊥

infix 4 _≢_
```

```
_≢_ : {A : Set} → A → A → Set
x ≢ y = ¬(x ≡ y)
```

The principle of *ex falso quodlibet* (EFQ) holds in Agda, in the sense that any type can be constructed from the bottom type. To show this, we do a case split on the instance of ⊥. There is no constructor for ⊥, which is stated using empty parentheses. Cases which are not constructable do not need proving.

```
⊥-elim : {A : Set} → ⊥ → A
⊥-elim ()
```

A proposition (type) is decidable if it can be proved (constructed), or otherwise if its proof (construction) leads to a proof (construction) of ⊥.

```
data Dec (A : Set) : Set where
  yes :   A → Dec A
  no  : ¬ A → Dec A
```

The constructors `yes` and `no` can be thought of as similar to the truth values `true` and `false` in the boolean type, with the addition that they keep the proof or disproof of the proposition for which they are acting as a truth value.

A unary predicate is *decidable* if each of its values is decidable.

```
Pred : Set → Set₁
Pred A = A → Set

Decidable : {A : Set} → Pred A → Set
Decidable P = ∀ x → Dec (P x)
```

The same could be defined for binary predicates, but this won't be needed. However, the special case of the equality predicate being decidable for a given type[1] will be used later.

```
Decidable≡ : Set → Set
Decidable≡ A = (x y : A) → Dec (x ≡ y)
```

Intuitively, inductively defined types which are not constructed from functions will have a decidable equality, simply by case analysis on the components from which they are constructed.

# 3 Nat.lagda

There is a built-in module for natural numbers, which defines the arithmetic operations and boolean relations, including a boolean-valued equality. We import and augment this with some propositions and predicates. The (unicode-renamed) definition of natural numbers is commented below.

```
open import Agda.Builtin.Nat renaming (Nat to ℕ) hiding (_<_) public

{-
  data ℕ : Set where
    zero : ℕ
    suc  : ℕ → ℕ
-}
```

[1]This is as much a property of the type as it is a property of the equality predicate for that type. A type with a decidable equality is called discrete [8]. In a programming-focused implementation, we could define a type expressing the discreteness of another type, and provide `instance` proofs for the discreteness of the types we work with. Agda would then find any required decidable equalities implicitly by instance resolution. However, in the interests of remaining closer to a type-theoretic implementation, we will eschew this Agda feature.

The built-in boolean-valued equality _==_ can be evaluated to check that `1 + 1 == 2` is `true`. However, this is not useful as a lemma. Instead, we would like to have a binary predicate for natural numbers which gives either a proof of equality or a proof of inequality. Such a predicate is itself a proof that equality of natural numbers is decidable, given the definition of `Decidable≡` above.

The proof is by case analysis on the arguments. In the case where both numbers are zero, they can be proven equal simply by `refl`. Where only one number is a successor, they can be proven not equal by doing case analysis on what their equality would be. As the only constructor for _≡_ requires that the left and right sides are the same, and `zero` cannot be unified with `suc _`, the cases are empty. Finally, if both numbers are successors, check if their predecessors are equal. If so, then equality follows. Otherwise, assuming the numbers are equal leads to a contradiction.

```
natEq : Decidable≡ ℕ
natEq zero    zero    = yes refl
natEq zero    (suc m) = no λ ()
natEq (suc n) zero    = no λ ()
natEq (suc n) (suc m) with natEq n m
...                   | yes refl = yes refl
...                   | no  n≢m  = no  λ { refl → n≢m refl }
```

A propositional order relation on the natural numbers can be defined as usual.

```
data _≤_ : ℕ → ℕ → Set where
  0≤n   : ∀{n} → zero ≤ n
  sn≤sm : ∀{n m} → n ≤ m → suc n ≤ suc m

_<_ : ℕ → ℕ → Set
n < m = suc n ≤ m
```

In the definition of '≤', the type is *indexed* by a pair of natural numbers, rather than parametrised (given specific names, on the left side of the colon). This is an example of a dependent type. The constructors do not produce values of the same type. Moreover, there are types for which there are no constructors. For example, there is no way of constructing `1 ≤ 0`. In this manner, dependent types can describe predicates.

The relation _≤_ is reflexive and transitive.

```
≤refl : ∀{n} → n ≤ n
≤refl {zero}  = 0≤n
≤refl {suc n} = sn≤sm ≤refl

≤trans : ∀{x y z} → x ≤ y → y ≤ z → x ≤ z
≤trans 0≤n         y≤z         = 0≤n
≤trans (sn≤sm x≤y) (sn≤sm y≤z) = sn≤sm (≤trans x≤y y≤z)
```

If $n < m$ then $m \nleq n$, and if $m \leq n$ then $n \nless m$. This can be expressed as a single proposition. To derive $\bot$, recurse on $n$ and $m$ until one of them is 0, at which point there is either no constructor for $n < m$ or no constructor for $m \leq n$.

```
ℕdisorder : ∀{n m} → n < m → m ≤ n → ⊥
ℕdisorder (sn≤sm n<m) (sn≤sm m≤n) = ℕdisorder n<m m≤n
```

Given natural numbers $n$ and $m$, it is possible to compute whether $n \leq m$ or $m \leq n$. To prove this, we first create a proposition `Compare n m` which is constructed by a proof of either of these.

```
data Compare (n m : ℕ) : Set where
  less : n ≤ m → Compare n m
  more : m ≤ n → Compare n m
```

It remains to show that given any $n$ and $m$, we may construct `Compare n m`.

```
compare : ∀ n m → Compare n m
compare zero    m       = less 0≤n
compare (suc n) zero    = more 0≤n
compare (suc n) (suc m) with compare n m
...                     | less n≤m = less (sn≤sm n≤m)
...                     | more m≤n = more (sn≤sm m≤n)
```

While it is possible to directly define a function which returns the greater of two natural numbers, this method preserves the proof showing which is greater. Defining a relation, and then supplying a function to construct it from all possible arguments is a common technique, and it will be used often.

# 4 List.lagda

We extend the built-in module for lists, by showing that if a predicate over a type is decidable, then given a list over that type, it is decidable if the predicate holds on any member, and it is decidable if the predicate holds on all members.

First, import the built-in list type. A simplified version of the definition is commented below.

```
open import Agda.Builtin.List public

{-
  data List (A : Set) : Set where
    []   : List A
    _::_ : A → List A → List A
-}
```

A list of type $A$ is either empty, or otherwise constructed by prepending an object of type $A$ to a list of type $A$. Given a predicate $P$ on $A$, the notion of $P$ holding on every element of a list can be defined in a similar way.

```
data All {A : Set} (P : Pred A) : List A → Set where
  []   : All P []
  _::_ : ∀{x xs} → P x → All P xs → All P (x :: xs)
```

In the case that $P$ is decidable, it is also decidable whether $P$ holds on every element of a list, by simply recursing through and examining $P$ on every element.

```
decAll : ∀{A P} → (p : Decidable P) → (xs : List A) → Dec (All P xs)
decAll p []        = yes []
decAll p (x :: xs) with p x
...                | no ¬Px = no λ { (Px :: _) → ¬Px Px }
...                | yes Px with decAll p xs
...                         | yes ∀xsP = yes (Px :: ∀xsP)
...                         | no ¬∀xsP = no  λ { (_ :: ∀xsP) → ¬∀xsP ∀xsP }
```

For $P$ to hold on *any* element of a list, it must either hold on the first element, or otherwise in the tail of the list.

```
data Any {A : Set} (P : Pred A) : List A → Set where
  [_]  : ∀{x xs} → P x → Any P (x :: xs)
  _::_ : ∀{xs} → (x : A) → Any P xs → Any P (x :: xs)
```

Again, the above is decidable for decidable predicates.

```
decAny : ∀{A P} → (p : Decidable P) → (xs : List A) → Dec (Any P xs)
decAny p []        = no λ ()
decAny p (x :: xs) with p x
...                | yes Px = yes [ Px ]
...                | no ¬Px with decAny p xs
...                          | yes ∃xsP = yes (x :: ∃xsP)
...                          | no ¬∃xsP = no  λ { [ Px ]       → ¬Px Px
                                                ; ( _ :: ∃xsP) → ¬∃xsP ∃xsP }
```

We can now define the membership predicate '∈' for lists; $x \in xs$ if any member of $xs$ is equal to $x$. The command `infix` sets the arity of the infix operators.

```
infix 4 _∈_ _∉_

_∈_ : {A : Set} → (x : A) → List A → Set
x ∈ xs = Any (x ≡_) xs

_∉_ : {A : Set} → (x : A) → List A → Set
x ∉ xs = ¬(x ∈ xs)
```

It follows that if equality is decidable, then membership is decidable.

```
dec∈ : ∀{A} → Decidable≡ A → (x : A) → (xs : List A) → Dec (x ∈ xs)
dec∈ _≟_ x xs = decAny (x ≟_) xs
```

# 5 Vec.lagda

Vectors are similar to lists, but the type is indexed by length. For example, vectors in $\mathbb{N}^2$ are of different type to vectors in $\mathbb{N}^3$.

```
data Vec (A : Set) : ℕ → Set where
  []   : Vec A zero
  _::_ : ∀{n} → A → Vec A n → Vec A (suc n)
```

We define All, Any, and membership the same was as for lists. The decidability proofs below are omitted, as they are identical to the corresponding proofs for lists.

```
data All {A : Set} (P : Pred A) : ∀{n} → Vec A n → Set where
  []   : All P []
  _::_ : ∀{x n} {xs : Vec A n} → P x → All P xs → All P (x :: xs)

decAll : ∀{A n P} → (p : Decidable P) → (xs : Vec A n) → Dec (All P xs)
-- Proof omitted.

data Any {A : Set} (P : Pred A) : ∀{n} → Vec A n → Set where
  [_]  : ∀{n x} {xs : Vec A n} → P x → Any P (x :: xs)
  _::_ : ∀{n}   {xs : Vec A n} → (x : A) → Any P xs → Any P (x :: xs)

decAny : ∀{A n P} → (p : Decidable P) → (xs : Vec A n) → Dec (Any P xs)
-- Proof omitted.
```

```
infix 4 _∈_ _∉_

_∈_ : {A : Set} {n : ℕ} → (x : A) → Vec A n → Set
x ∈ xs = Any (x ≡_) xs

_∉_ : {A : Set} {n : ℕ} → (x : A) → Vec A n → Set
x ∉ xs = ¬(x ∈ xs)

dec∈ : ∀{A n} → Decidable≡ A → (x : A) → (xs : Vec A n) → Dec (x ∈ xs)
dec∈ _≟_ x xs = decAny (x ≟_) xs
```

# 6 Formula.lagda

```
open import Agda.Builtin.Sigma
```

## 6.1 Basic definitions

We adopt the definitions of [10].

There are countably many variables, and there are countably many function symbols of each (natural) arity. Constants are functions with arity zero. Function symbols of different arities with the same index are considered distinct.

```
record Variable : Set where
  constructor var
  field
    varidx : ℕ

open Variable

record Function : Set where
  constructor func
  field
    funcidx   : ℕ
    funcarity : ℕ

open Function
```

Agda's `record` types fill the role of sigma types from MLTT.[2] By defining these as `record` types, we get destructors for accessing the indices and arities, which we then extract into the current module for ease of use. Note that the indices are natural numbers. While it seems equivalent and more natural to use string indices, strings are less useful for proofs. Internally, strings are not recursively defined as the natural numbers are; instead the string type is a postulated type which is bound to string literals.

Terms are either variables, or functions applied to the appropriate number of arguments (zero for constants).

```
data Term : Set where
  varterm  : Variable → Term
  functerm : (f : Function) → Vec Term (funcarity f) → Term
```

Relation symbols work the same way as function symbols.

[2] We will later define the type Σ for sigma types using `record` types.

```
record Relation : Set where
  constructor rel
  field
    idx   : ℕ
    arity : ℕ

open Relation renaming (idx to relidx ; arity to relarity)
```

A formula is either atomic (a prime formula), or formed from one of the logical connectives or quantifiers. We use 'Λ' (capital lambda) and 'V' (capital 'v') for '∀' and '∃', since '∀' is reserved by Agda.[3]

```
data Formula : Set where
  atom    : (r : Relation) → Vec Term (relarity r) → Formula
  _⇒_     : Formula  → Formula → Formula
  _∧_     : Formula  → Formula → Formula
  _∨_     : Formula  → Formula → Formula
  Λ       : Variable → Formula → Formula
  V       : Variable → Formula → Formula

_⇔_ : Formula → Formula → Formula
Φ ⇔ Ψ = (Φ ⇒ Ψ) ∧ (Ψ ⇒ Φ)
```

The logical connectives are right-associative, and have the usual order of precedence.

```
infixr 105 _⇒_ _⇔_
infixr 106 _∨_
infixr 107 _∧_
```

Equality of formulae is decidable. Logically, this follows from the fact that formulae are inductively defined. The proof is obtained by case analysis, using lemmas on the types used to construct formulae. As these proofs are unremarkable, and follow the same pattern as the proof for decidable equality of natural numbers above, they are omitted.

```
varEq : Decidable≡ Variable
-- Proof omitted.

relEq : Decidable≡ Relation
-- Proof omitted.

funcEq : Decidable≡ Function
-- Proof omitted.

termEq : Decidable≡ Term
-- Proof omitted.

formulaEq : Decidable≡ Formula
-- Proof omitted.
```

## 6.2 Variable freedom

We define the conditions for a variable to be *not free* in a formula. Instead of first defining *free* and then taking *not free* to be the negation, we use a positive definition for not free, since later definitions only ever require proof that a variable is not free.

[3] While the typical n-ary logical operator symbols '⋁' and '⋀' are available, they are more easily confused with the symbols '∧' and '∨' for 'and' and 'or', and are unavailable in some fonts.

For a given term $t$, $x$ is not in $t$ if $t$ is a variable other than $x$. Otherwise if the term is a function on arguments $ts$, then $x$ is not in $t$ if it is not anywhere in $ts$, which can be checked by applying `All` to this definition. Separating the declaration and definition of `_NotInTerm_` allows it to be defined mutually with the case for a vector of terms.

```
data _NotInTerm_ (x : Variable) : Term → Set

_NotInTerms_ : ∀{n} → Variable → Vec Term n → Set
x NotInTerms ts = All (x NotInTerm_) ts

data _NotInTerm_ x where
  varterm  : ∀{y} → x ≢ y → x NotInTerm (varterm y)
  functerm : ∀{f} {us : Vec Term (funcarity f)}
             → x NotInTerms us → x NotInTerm (functerm f us)
```

A variable is now not free in a formula according to the obvious recursive definition. It is not free inside an atom if it is not inside that atom, meaning it is not in the terms that the relation is operating on. It is not free inside a quantification over a subformula either if it is the quantification variable, or else if it is not free in the subformula. Separate constructors are given for each case.

```
data _NotFreeIn_ : Variable → Formula → Set where
  atom : ∀{x r} {ts : Vec Term (relarity r)}
             → x NotInTerms ts → x NotFreeIn (atom r ts)
  _⇒_  : ∀{x α β} → x NotFreeIn α → x NotFreeIn β → x NotFreeIn (α ⇒ β)
  _∧_  : ∀{x α β} → x NotFreeIn α → x NotFreeIn β → x NotFreeIn (α ∧ β)
  _∨_  : ∀{x α β} → x NotFreeIn α → x NotFreeIn β → x NotFreeIn (α ∨ β)
  Λ↓   : ∀ x α     → x NotFreeIn Λ x α
  V↓   : ∀ x α     → x NotFreeIn V x α
  Λ    : ∀{x α}    → ∀ y → x NotFreeIn α → x NotFreeIn Λ y α
  V    : ∀{x α}    → ∀ y → x NotFreeIn α → x NotFreeIn V y α
```

**Lemma 6.2.1.** *Variable occurrence within a vector of terms is decidable.*

*Proof.* Search through the vector for occurrences of the variable. In the following code we will use names like `x∉t` to denote proofs of '$x$ is not in term $t$', `x∉ts` for '$x$ is not in any terms in $ts$', and `x∉α` for '$x$ is not free in $\alpha$'.

```
_notInTerms_ : ∀{n} → ∀ x → (ts : Vec Term n) → Dec (x NotInTerms ts)
x notInTerms [] = yes []
```

To check against a variable term, use the decidable equality of variables, then recurse over the rest of the terms.

```
x notInTerms (varterm y :: ts) with varEq x y
...     | yes refl = no λ { (varterm x≢x :: _) → x≢x refl }
...     | no  x≢y  with x notInTerms ts
...                | yes x∉ts = yes (varterm x≢y :: x∉ts)
...                | no ¬x∉ts = no  λ { (_ :: x∉ts) → ¬x∉ts x∉ts }
```

To check against a function term, recurse over the arguments, then recurse over the rest of the terms.

```
x notInTerms (functerm f us :: ts) with x notInTerms us
...     | no ¬x∉us = no λ { (functerm x∉us :: _) → ¬x∉us x∉us }
...     | yes x∉us with x notInTerms ts
...                | yes x∉ts = yes (functerm x∉us :: x∉ts)
...                | no ¬x∉ts = no  λ { (_ :: x∉ts) → ¬x∉ts x∉ts }
```

Each case checks if $x$ is free in the remaining terms in the vector. A shorter proof would do this check at the same time as doing a case split on the first term. However, if a term for which $x$ is free is found, it is not necessary to continue recursing through the vector, so it is better computationally not to do so. □

The same logic can be used for a single term, calling the above function to check function arguments. The proposition `_NotInTerms_` is defined using `All` and `_NotInTerm_`, so it is tempting to try to first prove that the single term case is decidable, and then generalise to vectors using the lemma that `All` is decidable for decidable predicates. However, this would not be structurally recursive, and so Agda would not see this as terminating. Above, the case `x notInTerms t :: ts` depends on the result of `x notInTerms ts`, which is in fact primitively recursive. However, if it instead depended on the result of `all (x notInTerm_) ts`, Agda cannot determine that `x notInTerm_` will be applied only to arguments structurally smaller than `t :: ts`.

```
_notInTerm_ : (x : Variable) → (t : Term) → Dec (x NotInTerm t)
x notInTerm varterm y     with varEq x y
...                         | yes refl = no  λ { (varterm x≢x) → x≢x refl }
...                         | no  x≢y  = yes (varterm x≢y)
x notInTerm functerm f ts with x notInTerms ts
...                         | yes x∉ts = yes (functerm x∉ts)
...                         | no ¬x∉ts = no  λ { (functerm x∉ts) → ¬x∉ts x∉ts }
```

**Proposition 6.2.2.** *Variable freedom is decidable.*

*Proof.* For atoms, apply the lemma above. Otherwise, check recursively, checking if the variable matches the quantifying variable in the case of quantifiers.

```
_notFreeIn_ : (x : Variable) → (α : Formula) → Dec (x NotFreeIn α)
x notFreeIn atom r ts with x notInTerms ts
...                     | yes x∉ts = yes (atom x∉ts)
...                     | no ¬x∉ts = no  λ { (atom x∉ts) → ¬x∉ts x∉ts }
x notFreeIn (α ⇒ β)    with x notFreeIn α | x notFreeIn β
...                     | yes x∉α | yes x∉β = yes (x∉α ⇒ x∉β)
...                     | no ¬x∉α | _       = no  λ { (x∉α ⇒ _  ) → ¬x∉α x∉α }
...                     | _       | no ¬x∉β = no  λ { (_   ⇒ x∉β) → ¬x∉β x∉β }
x notFreeIn (α ∧ β)    with x notFreeIn α | x notFreeIn β
...                     | yes x∉α | yes x∉β = yes (x∉α ∧ x∉β)
...                     | no ¬x∉α | _       = no  λ { (x∉α ∧ _  ) → ¬x∉α x∉α }
...                     | _       | no ¬x∉β = no  λ { (_   ∧ x∉β) → ¬x∉β x∉β }
x notFreeIn (α ∨ β)    with x notFreeIn α | x notFreeIn β
...                     | yes x∉α | yes x∉β = yes (x∉α ∨ x∉β)
...                     | no ¬x∉α | _       = no  λ { (x∉α ∨ _  ) → ¬x∉α x∉α }
...                     | _       | no ¬x∉β = no  λ { (_   ∨ x∉β) → ¬x∉β x∉β }
x notFreeIn Λ  y α     with varEq x y
...                     | yes refl = yes (Λ↓ x α)
...                     | no  x≢y  with x notFreeIn α
...                                  | yes x∉α = yes (Λ y x∉α)
...                                  | no ¬x∉α = no  λ { (Λ↓ x α)   → x≢y refl
                                                       ; (Λ y x∉α) → ¬x∉α x∉α }
x notFreeIn V  y α     with varEq x y
...                     | yes refl = yes (V↓ x α)
...                     | no  x≢y  with x notFreeIn α
...                                  | yes x∉α = yes (V y x∉α)
...                                  | no ¬x∉α = no  λ { (V↓ x α)   → x≢y refl
                                                       ; (V y x∉α) → ¬x∉α x∉α }
```
□

## 6.3 Substitutions

We define what it means for a formula $\beta$ to be obtained from $\alpha$ by replacing all free instances of a variable $x$ with term $t$, by giving a relation `_[_/_]≡_`. Some of the natural deduction rules will involve variable substitution, and the type of the result of the deduction will depend on the result of the substitution. If we instead defined substitution as a function instead of a relation, we would have to provide equality proofs about the value computed by the function for (sometimes arbitrary) formulae. This is unweildy, and cannot be solved in general by Agda's proof search. Instead, we will define our relation so that it can be proved easily (and automatically) when doing natural deduction, and then later give a function which computes both a formula $\beta$, and a proof that $\beta$ is the required substitution.

The definitions below have a similar structure to that of `_NotFreeIn_` above. The more general case of replacing terms with terms is not needed for natural deduction.

Inside a vector of terms, wherever $x$ occurs, it is replaced with $t$. Any variable distinct from $x$ is left unchanged. For a function term, $x$ is replaced with $t$ inside all of the arguments.

```
data [_][_/_]≡_ : ∀{n} → Vec Term n → Variable → Term → Vec Term n → Set

data ⟨_⟩[_/_]≡_ : Term → Variable → Term → Term → Set where
  varterm≡ : ∀{x t}          → ⟨ varterm x ⟩[ x / t ]≡ t
  varterm≢ : ∀{x t y}        → x ≢ y → ⟨ varterm y ⟩[ x / t ]≡ varterm y
  functerm : ∀{x t f us vs} → [ us ][ x / t ]≡ vs
                              → ⟨ functerm f us ⟩[ x / t ]≡ functerm f vs

data [_][_/_]≡_ where
  []   : ∀{x t} → [ [] ][ x / t ]≡ []
  _∷_ : ∀{x t u v n} {us vs : Vec Term n}
        → ⟨ u ⟩[ x / t ]≡ v → [ us ][ x / t ]≡ vs
        → [ u ∷ us ][ x / t ]≡ (v ∷ vs)
```

The definition for formulae follows.

```
data _[_/_]≡_ : Formula → Variable → Term → Formula → Set where
```

The `ident` constructor gives the case that replacing $x$ with $x$ yields the original formula. While this can be proved as a derived rule, in practice it is the case we usually want to use. Providing a constructor allows Agda's proof search to apply this case easily.

```
  ident    : ∀ α x → α [ x / varterm x ]≡ α
```

If $x$ is not free in $\alpha$, then replacing it with any term should leave $\alpha$ unchanged. This rule is not derivable when $t$ is not otherwise able to be substituted for $x$ in $\alpha$. For example, without this constructor it would not be possible to prove that $(\forall y A)[x/y] \equiv (\forall y A)$, where $A$ is a propositional formula.

```
  notfree : ∀{α x t} → x NotFreeIn α → α [ x / t ]≡ α
```

The propositional cases are similar to those of the `_NotFreeIn_` type above.

```
  atom     : ∀{x t}
               → (r : Relation) → {xs ys : Vec Term (relarity r)}
               → [ xs ][ x / t ]≡ ys → (atom r xs) [ x / t ]≡ (atom r ys)
  _⇒_      : ∀{α α′ β β′ x t}
               → α [ x / t ]≡ α′ → β [ x / t ]≡ β′
               → (α ⇒ β) [ x / t ]≡ (α′ ⇒ β′)
  _∧_      : ∀{α α′ β β′ x t}
               → α [ x / t ]≡ α′ → β [ x / t ]≡ β′
```

```
              → (α ∧ β) [ x / t ]≡ (α′ ∧ β′)
_∨_       : ∀{α α′ β β′ x t}
              → α [ x / t ]≡ α′ → β [ x / t ]≡ β′
              → (α ∨ β) [ x / t ]≡ (α′ ∨ β′)
```

Variable substitution for a quantified formula has two cases, which are similar to their counterparts in _NotFreeIn_. If $x$ is the quantification variable, then the formula is unchanged.

```
Λ↓        : ∀{t} → ∀ x α → (Λ x α) [ x / t ]≡ (Λ x α)
V↓        : ∀{t} → ∀ x α → (V x α) [ x / t ]≡ (V x α)
```

Finally, if $x$ is not the quantification variable, and the quantification variable does not appear in $t$, then the substitution simply occurs inside the quantification.

```
Λ         : ∀{α β x y t} → x ≢ y → y NotInTerm t
              → α [ x / t ]≡ β → (Λ y α) [ x / t ]≡ (Λ y β)
V         : ∀{α β x y t} → x ≢ y → y NotInTerm t
              → α [ x / t ]≡ β → (V y α) [ x / t ]≡ (V y β)
```

Given $\alpha$, $x$, $t$, the $\beta$ satisfying $\alpha[x/t] \equiv \beta$ should be unique, so that variable substitution is functional. This can first be shown for the special cases `ident` and `notfree`, by recursing through the constructors down to the atomic case, and recursing through the term substitutions down to the variable terms. The proofs simply have `refl` on the right side of every line, and are omitted. Their structures are very similar to the two proofs that follow afterward.

```
subIdentFunc : ∀{α x β} → α [ x / varterm x ]≡ β → α ≡ β
-- Proof omitted.
```

```
subNotFreeFunc : ∀{α x t β} → α [ x / t ]≡ β → x NotFreeIn α → α ≡ β
-- Proof omitted.
```

**Lemma 6.3.1.** *Variable substitution inside a vector of terms is functional.*

*Proof.* The constructors for term substitution have no overlap.

```
subTermsFunc : ∀{n x t} {us vs ws : Vec Term n}
               → [ us ][ x / t ]≡ vs → [ us ][ x / t ]≡ ws → vs ≡ ws
subTermsFunc [] [] = refl
```

First recurse over the rest of the two vectors.

```
subTermsFunc (s :: ss) (r :: rs) with subTermsFunc ss rs
```

It is possible to pattern match inside the `with` block to examine the two substitutions made for the heads of the vectors. In the case that the first term is substituted using `varterm≡` in each case, the resulting vectors must both have $x$ at the head, so the proof is `refl`.

```
subTermsFunc (varterm≡     :: _) (varterm≡     :: _) | refl = refl
```

It would be contradictory for the first term in $us$ to both match and differ from $x$.

```
subTermsFunc (varterm≡     :: _) (varterm≢ x≢x :: _) | refl = ⊥-elim (x≢x refl)
subTermsFunc (varterm≢ x≢x :: _) (varterm≡     :: _) | refl = ⊥-elim (x≢x refl)
```

If the head of $us$ is a variable different from $x$, then it is unchanged in each case, so the proof is `refl`.

```
subTermsFunc (varterm≢ x≢y :: _) (varterm≢ _   :: _) | refl = refl
```

Finally, in the case of a function, recurse over the vector of arguments. The `rewrite` construction uses a proof of equality to unify terms. It is an abbreviation for doing `with`-abstraction on a proof of `refl`.

```
subTermsFunc (functerm st  :: _) (functerm rt  :: _)
             | refl rewrite subTermsFunc st rt = refl
```
□

**Proposition 6.3.2.** *Variable substitution is functional.*

*Proof.*

```
subFunc : ∀{x t α β γ} → α [ x / t ]≡ β → α [ x / t ]≡ γ → β ≡ γ
```

If either substitution came from `ident` or `notfree`, invoke one of the above lemmas. If they occurred in the right substitution, the lemmas prove $\gamma \equiv \beta$, so `rewrite` is used to recover $\beta \equiv \gamma$.

```
subFunc (ident α x)   s             = subIdentFunc s
subFunc (notfree x∉α) s             = subNotFreeFunc s x∉α
subFunc r             (ident α x)   rewrite subIdentFunc r       = refl
subFunc r             (notfree x∉α) rewrite subNotFreeFunc r x∉α = refl
```

The atomic case comes from the previous lemma.

```
subFunc (atom p r)    (atom .p s)   rewrite subTermsFunc r s = refl
```

The propositional connectives can be proved inductively.

```
subFunc (r₁ ⇒ r₂)     (s₁ ⇒ s₂)     with subFunc r₁ s₁ | subFunc r₂ s₂
...                                 | refl | refl = refl
subFunc (r₁ ∧ r₂)     (s₁ ∧ s₂)     with subFunc r₁ s₁ | subFunc r₂ s₂
...                                 | refl | refl = refl
subFunc (r₁ ∨ r₂)     (s₁ ∨ s₂)     with subFunc r₁ s₁ | subFunc r₂ s₂
...                                 | refl | refl = refl
```

If the formula is a quantification over $x$, then neither substitution changes the formula.

```
subFunc (Λ↓ x α)      (Λ↓ .x .α)    = refl
subFunc (V↓ x α)      (V↓ .x .α)    = refl
```

It is contradictory for one substitution to occur by matching $x$ with the quantifier variable, and the other to have a different quantifier.

```
subFunc (Λ↓ x α)      (Λ x≢x _ s)   = ⊥-elim (x≢x refl)
subFunc (V↓ x α)      (V x≢x _ s)   = ⊥-elim (x≢x refl)
subFunc (Λ x≢x _ r)   (Λ↓ x α)      = ⊥-elim (x≢x refl)
subFunc (V x≢x _ r)   (V↓ x α)      = ⊥-elim (x≢x refl)
```

Finally, if the formula is a quantification over a variable other than $x$, then substitution occurs inside the quantified formula, so recurse inside those substitutions.

```
subFunc (Λ _ _ r)     (Λ _ _ s)     rewrite subFunc r s = refl
subFunc (V _ _ r)     (V _ _ s)     rewrite subFunc r s = refl
```
□

We have now shown that substitution is functional, and so would like to construct a function that computes substitutions. However, substitutions do not always exist. For example, there is no way of constructing a formula for $(\forall y Px)[x/y]$. In general, $\alpha[x/t]$ exists only if $t$ is *free for $x$ in $\alpha$*, meaning no variables in $t$ would become bound inside $\alpha$. This can be formalised by using (with minor modification) the rules of [14].

```
data _FreeFor_In_ (t : Term) (x : Variable) : Formula → Set where
  notfree : ∀{α} → x NotFreeIn α → t FreeFor x In α
  atom    : ∀ r us → t FreeFor x In atom r us
  _⇒_     : ∀{α β} → t FreeFor x In α → t FreeFor x In β
                   → t FreeFor x In α ⇒ β
  _∧_     : ∀{α β} → t FreeFor x In α → t FreeFor x In β
                   → t FreeFor x In α ∧ β
  _∨_     : ∀{α β} → t FreeFor x In α → t FreeFor x In β
                   → t FreeFor x In α ∨ β
  Λ↓      : ∀ α → t FreeFor x In Λ x α
  V↓      : ∀ α → t FreeFor x In V x α
  Λ       : ∀{α y} → y NotInTerm t → t FreeFor x In α → t FreeFor x In Λ y α
  V       : ∀{α y} → y NotInTerm t → t FreeFor x In α → t FreeFor x In V y α
```

The definitions above for variable substitution lead directly to a procedure for computing substitutions. Given $\alpha$, $x$, $t$, and a proof that $t$ is free for $x$ in $\alpha$, we compute a $\beta$ and a proof that $\alpha[x/t] \equiv \beta$.

The built-in sigma (dependent sum) type has been imported. A simplified version of its definition is commented below.

```
{-
  record Σ (A : Set) (B : A → Set) : Set where
    constructor _,_
    field
      fst : A
      snd : B fst
-}
```

A proof of a sigma type encapsulates both a value and a proof regarding that value. Proposition $\Sigma A B$ can be proved by providing an $x$ of type $A$, and a proof of $Bx$. This means that the sigma type can be used to define existential propositions.

**Lemma 6.3.3.** *Every vector of terms has a substitution of any variable with any term.*

*Proof.* Recurse through all function arguments, and replace any variables equal to $x$ with $t$. We do a case split on the first term, and use a `with` block to get the substitution for the rest of the vector simultaneously, since this substitution is required in either case.

```
[_][_/_] : ∀{n} → (us : Vec Term n) → ∀ x t → Σ _ [ us ][ x / t ]≡_
[ []                    ][ x / t ] = [] , []
[ u               :: us ][ x / t ] with [ us ][ x / t ]
[ varterm y       :: us ][ x / t ] | vs , vspf with varEq x y
...       | yes refl   = (t              :: vs) , (varterm≡      :: vspf)
...       | no  x≢y    = (varterm y      :: vs) , (varterm≢ x≢y  :: vspf)
[ functerm f ws :: us ][ x / t ] | vs , vspf with [ ws ][ x / t ]
...       | xs , xspf = (functerm f xs :: vs) , (functerm xspf :: vspf)
```

□

**Proposition 6.3.4.** *If $t$ is free for $x$ in $\alpha$, then there is a substitution of $x$ with $t$ in $\alpha$.*

*Proof.* The proof that $t$ is free for $x$ in formula must be supplied. The term $t$ is fixed by supplying such a proof, so for convenience of notation, the proof is supplied in place of the term.

```
_[_/_] : ∀{t} → ∀ α x → t FreeFor x In α → Σ Formula (α [ x / t ]≡_)
α [ x / notfree ¬x∉α ]          = α , notfree ¬x∉α
```

For atomic formulae, apply the above lemma.

```
_[_/_] {t} (atom r ts) x tff with [ ts ][ x / t ]
...                           | ts' , tspf = atom r ts' , atom r tspf
```

For the propositional connectives, the substitution is obtained recursively.

```
(α ⇒ β) [ x / tffα ⇒ tffβ ]  with α [ x / tffα ] | β [ x / tffβ ]
...                           | α' , αpf | β' , βpf = α' ⇒ β' , αpf ⇒ βpf
(α ∧ β) [ x / tffα ∧ tffβ ]  with α [ x / tffα ] | β [ x / tffβ ]
...                           | α' , αpf | β' , βpf = α' ∧ β' , αpf ∧ βpf
(α ∨ β) [ x / tffα ∨ tffβ ]  with α [ x / tffα ] | β [ x / tffβ ]
...                           | α' , αpf | β' , βpf = α' ∨ β' , αpf ∨ βpf
```

For generalisation, check if $x$ is the quantifier variable, and if so do nothing. Otherwise, recurse.

```
Λ y α [ .y / Λ↓ .α ]          = Λ y α , Λ↓ y α
V y α [ .y / V↓ .α ]          = V y α , V↓ y α
Λ y α [ x / Λ y∉t tffα ]      with varEq x y
...                           | yes refl = Λ y α , Λ↓ y α
...                           | no  x≢y  with α [ x / tffα ]
...                                      | α' , αpf = Λ y α' , Λ x≢y y∉t αpf
V y α [ x / V y∉t tffα ]      with varEq x y
...                           | yes refl = V y α , V↓ y α
...                           | no  x≢y  with α [ x / tffα ]
...                                      | α' , αpf = V y α' , V x≢y y∉t αpf
```
□

We have proved that if $t$ is free for $x$ in $\alpha$ then $\alpha[x/t]$ exists. The converse is also true, meaning that `_FreeFor_In_` precisely captures the notion of a substitution being possible. The proof is straightforward by induction on formula substitution, with the base case of atomic formulae being trivial.

```
subFreeFor : ∀{α x t β} → α [ x / t ]≡ β → t FreeFor x In α
-- Proof omitted.
```

**Proposition 6.3.5.** *If a variable has been substituted by a term not involving that variable, then the variable is not free in the resulting formula.*

*Proof.*

```
subNotFree : ∀{α x t β} → x NotInTerm t → α [ x / t ]≡ β → x NotFreeIn β
```

The case where the substitution was constructed by `ident` is absurd, since $x$ can't not be in term $x$.

```
subNotFree (varterm x≢x) (ident α x)    = ⊥-elim (x≢x refl)
```

If the substitution was constructed by `notfree`, then $\alpha = \beta$, so $x$ is not free in $\beta$.

```
subNotFree x∉t            (notfree x∉α)  = x∉α
```

For atomic formulae, we use an inline lemma that the proposition holds for vectors of terms. Every variable in a term is either equal to $x$, and so gets replaced with $t$, or else differs from $x$.

```
subNotFree x∉t (atom r subts)   = atom (φ x∉t subts)
  where
    φ : ∀{n x t} {us vs : Vec Term n}
        → x NotInTerm t → [ us ][ x / t ]≡ vs → x NotInTerms vs
    φ x∉t []                        = []
    φ x∉t (varterm≡      :: subus)  = x∉t                    :: φ x∉t subus
    φ x∉t (varterm≢ neq :: subus)  = varterm neq            :: φ x∉t subus
    φ x∉t (functerm sub :: subus)  = functerm (φ x∉t sub) :: φ x∉t subus
```

The remaining cases follow by recursion.

```
subNotFree x∉t (subα ⇒ subβ)   = subNotFree x∉t subα ⇒ subNotFree x∉t subβ
subNotFree x∉t (subα ∧ subβ)   = subNotFree x∉t subα ∧ subNotFree x∉t subβ
subNotFree x∉t (subα ∨ subβ)   = subNotFree x∉t subα ∨ subNotFree x∉t subβ
subNotFree x∉t (Λ↓ y α)        = Λ↓ y α
subNotFree x∉t (Λ x≢y y∉t sub) = Λ _ (subNotFree x∉t sub)
subNotFree x∉t (V↓ y α)        = V↓ y α
subNotFree x∉t (V x≢y y∉t sub) = V _ (subNotFree x∉t sub)
```

□

**Proposition 6.3.6.** *Substituting with a variable which is not free is invertible by reversing the substitution.*

*Proof.*

```
subInverse : ∀{ω α x β} → ω NotFreeIn α
                        → α [ x / varterm ω ]≡ β → β [ ω / varterm x ]≡ α
```

The cases where the substitution was obtained with the `ident` or `notfree` constructors are trivial, since the formula has not been changed.

```
subInverse _           (ident α x)     = ident α x
subInverse ω∉α         (notfree x∉α)   = notfree ω∉α
```

In the atomic case, we use an inline lemma that the proposition holds for vectors of terms.

```
subInverse (atom x∉ts) (atom r subts) = atom r (ψ x∉ts subts)
  where
    ψ : ∀{n x ω} {us vs : Vec Term n}
        → ω NotInTerms us → [ us ][ x / varterm ω ]≡ vs
        → [ vs ][ ω / varterm x ]≡ us
    ψ ω∉us                     []                       = []
    ψ (_            :: ω∉us) (varterm≡      :: subus) = varterm≡
                                                         :: ψ ω∉us subus
    ψ (varterm ω≢y  :: ω∉us) (varterm≢ x≢ω :: subus) = varterm≢ ω≢y
                                                         :: ψ ω∉us subus
    ψ (functerm ω∉ts :: ω∉us) (functerm sub :: subus) = functerm (ψ ω∉ts sub)
                                                         :: ψ ω∉us subus
```

The propositional connective cases are solved by recursion.

```
subInverse (ω∉α ⇒ ω∉β) (sα ⇒ sβ) = subInverse ω∉α sα ⇒ subInverse ω∉β sβ
subInverse (ω∉α ∧ ω∉β) (sα ∧ sβ) = subInverse ω∉α sα ∧ subInverse ω∉β sβ
subInverse (ω∉α ∨ ω∉β) (sα ∨ sβ) = subInverse ω∉α sα ∨ subInverse ω∉β sβ
```

If the substitution changed nothing because the substitution variable was a quantifier variable, then $\omega$ is still not free in $\beta$.

```
subInverse ω∉α (Λ↓ x α) = notfree ω∉α
subInverse ω∉α (V↓ x α) = notfree ω∉α
```

Now consider the case where the substitution occurred inside a quantifier. It is absurd for $\omega$ to be the quantifer, since it would not have been allowed to substitute $x$ with $\omega$.

```
subInverse (Λ↓ x α) (Λ _ (varterm x≢x) _) = ⊥-elim (x≢x refl)
subInverse (V↓ x α) (V _ (varterm x≢x) _) = ⊥-elim (x≢x refl)
```

Suppose the formula was $\forall y \alpha$. Again discard the case where $\omega$ is $y$.

```
subInverse {ω} (Λ y ω∉α) (Λ _ y∉ω              _) with varEq ω y
subInverse {ω} (Λ y ω∉α) (Λ _ (varterm y≢y) _) | yes refl = ⊥-elim (y≢y refl)
```

Recurse inside the quantifier, turning a proof of $x \neq y$ into $y \neq x$.

```
subInverse {ω} (Λ y ω∉α) (Λ x≢y y∉ω sub)       | no  ω≢y
    = Λ ω≢y (varterm λ { refl → x≢y refl }) (subInverse ω∉α sub)
```

The same applies if the formula was $\exists y\alpha$.

```
subInverse {ω} (V y ω∉α) (V _ y∉ω              _) with varEq ω y
subInverse {ω} (V y ω∉α) (V _ (varterm y≢y) _) | yes refl = ⊥-elim (y≢y refl)
subInverse {ω} (V y ω∉α) (V x≢y y∉ω sub)       | no  ω≢y
    = V ω≢y (varterm λ { refl → x≢y refl }) (subInverse ω∉α sub)
```

□

## 6.4 Fresh variables

A variable is *fresh* if appears nowhere (free or bound) in a formula.

```
data _FreshIn_ (x : Variable) : Formula → Set where
  atom : ∀{r ts} → x NotInTerms ts → x FreshIn (atom r ts)
  _⇒_  : ∀{α β}  → x FreshIn α → x FreshIn β → x FreshIn α ⇒ β
  _∧_  : ∀{α β}  → x FreshIn α → x FreshIn β → x FreshIn α ∧ β
  _∨_  : ∀{α β}  → x FreshIn α → x FreshIn β → x FreshIn α ∨ β
  Λ    : ∀{α y}  → y ≢ x → x FreshIn α → x FreshIn Λ y α
  V    : ∀{α y}  → y ≢ x → x FreshIn α → x FreshIn V y α
```

Certainly, if a variable is fresh in a formula, then it is also not free, and every term is free for that variable. The proofs are trivial, and are omitted.

```
freshNotFree : ∀{α x} → x FreshIn α → x NotFreeIn α
-- Proof omitted.

freshFreeFor : ∀{α x} → x FreshIn α → ∀ y → (varterm x) FreeFor y In α
-- Proof omitted.
```

For the purposes of variable substitution, we will later need a way to generate a fresh variable for a given formula. Only finitely many variables occur in a given term or formula, so there is a greatest (with respect to the natural number indexing) variable occurring in each term or formula; all variables greater than this are fresh. We will first compute this variable, and then use its successor as the fresh variable. This means that the least fresh variable will not be found. For example, for $Px_0 \lor Px_2$, we find that $x_3, x_4, \ldots$ are fresh, missing $x_1$. However, finding the least fresh variable cannot be done with a simple recursive procedure. Consider $\alpha = (Px_0 \lor Px_2) \land Px_1$; we find $x_1$ is fresh to the left of the conjunctive, and $x_0$ is fresh to the right, but this does not indicate that $x_2$ will not be fresh in $\alpha$.

**Lemma 6.4.1.** *There is an upper bound on the variables occurring in a given vector of terms.*

*Proof.* We call this function `maxVarTerms`, but will not actually prove that this is the least upper bound in particular.

```
maxVarTerms : ∀{k} → (ts : Vec Term k)
              → Σ Variable (λ ⌜ts⌝
                           → ∀ n → varidx ⌜ts⌝ < n → var n NotInTerms ts)
maxVarTerms []                 = var zero , (λ _ _ → [])
```

If the first term is a variable, check if its index is greater than or equal to the greatest variable in the rest of the terms.

```
maxVarTerms (varterm x     :: ts) with maxVarTerms ts
... | ⌜ts⌝ , tspf with compare (varidx x) (varidx ⌜ts⌝)
```

If so, use it.

```
...                 | more ⌜ts⌝≤x = x , maxIsx
  where
    maxIsx : ∀ n → (varidx x) < n → (var n) NotInTerms (varterm x :: ts)
    maxIsx n x<n    = varterm (λ { refl → ℕdisorder x<n ≤refl })
                      :: tspf n (≤trans (sn≤sm ⌜ts⌝≤x) x<n)
```

Otherwise, use the greatest variable in the rest of the terms.

```
...                 | less x≤⌜ts⌝ = ⌜ts⌝ , ⌜ts⌝pf
  where
    ⌜ts⌝pf : ∀ n → varidx ⌜ts⌝ < n → var n NotInTerms (varterm x :: ts)
    ⌜ts⌝pf n ⌜ts⌝<n = varterm (λ { refl → ℕdisorder ⌜ts⌝<n x≤⌜ts⌝ })
                      :: tspf n ⌜ts⌝<n
```

If the first term is a function, then check if the greatest variable in its arguments is greater than or equal to the greatest variable of the rest of the terms.

```
maxVarTerms (functerm f us :: ts) with maxVarTerms us | maxVarTerms ts
... | ⌜us⌝ , uspf | ⌜ts⌝ , tspf with compare (varidx ⌜us⌝) (varidx ⌜ts⌝)
```

If so, use it.

```
...                            | more ⌜ts⌝≤⌜us⌝ = ⌜us⌝ , ⌜us⌝pf
  where
    ⌜us⌝pf : ∀ n → varidx ⌜us⌝ < n → (var n) NotInTerms (functerm f us :: ts)
    ⌜us⌝pf n ⌜us⌝<n = functerm (uspf n ⌜us⌝<n)
                      :: tspf n (≤trans (sn≤sm ⌜ts⌝≤⌜us⌝) ⌜us⌝<n)
```

If not, use the greatest variable in the rest of the terms.

```
...                            | less ⌜us⌝≤⌜ts⌝ = ⌜ts⌝ , ⌜ts⌝pf
  where
    ⌜ts⌝pf : ∀ n → varidx ⌜ts⌝ < n → (var n) NotInTerms (functerm f us :: ts)
    ⌜ts⌝pf n ⌜ts⌝<n = functerm (uspf n (≤trans (sn≤sm ⌜us⌝≤⌜ts⌝) ⌜ts⌝<n))
                      :: tspf n ⌜ts⌝<n
```
□

**Proposition 6.4.2.** *There is an upper bound on the variables occurring in a given formula.*

*Proof.*

```
maxVar : ∀ α → Σ Variable λ ⌜α⌝ → ∀ n → varidx ⌜α⌝ < n → var n FreshIn α
```

In the atomic case, apply the above lemma to find the greatest variable occuring.

```
maxVar (atom r ts) with maxVarTerms ts
...                  | ⌜ts⌝ , tspf = ⌜ts⌝ , λ n ⌜ts⌝<n → atom (tspf n ⌜ts⌝<n)
```

If all variables greater than $\lceil\alpha\rceil$ are fresh in $\alpha$, and all greater than $\lceil\beta\rceil$ are fresh in $\beta$, then any variable greater than $\max\{\lceil\alpha\rceil, \lceil\beta\rceil\}$ will be fresh in $\alpha \rightarrow \beta$.

```
maxVar (α ⇒ β) with maxVar α | maxVar β
...     | ⌜α⌝ , αpf | ⌜β⌝ , βpf with compare (varidx ⌜α⌝) (varidx ⌜β⌝)
...                              | less ⌜α⌝≤⌜β⌝ = ⌜β⌝ , maxIs⌜β⌝
    where
      maxIs⌜β⌝ : ∀ n → varidx ⌜β⌝ < n → var n FreshIn (α ⇒ β)
      maxIs⌜β⌝ n ⌜β⌝<n = αpf n (≤trans (sn≤sm ⌜α⌝≤⌜β⌝) ⌜β⌝<n) ⇒ βpf n ⌜β⌝<n
...                              | more ⌜β⌝≤⌜α⌝ = ⌜α⌝ , maxIs⌜α⌝
    where
      maxIs⌜α⌝ : ∀ n → varidx ⌜α⌝ < n → var n FreshIn (α ⇒ β)
      maxIs⌜α⌝ n ⌜α⌝<n = αpf n ⌜α⌝<n ⇒ βpf n (≤trans (sn≤sm ⌜β⌝≤⌜α⌝) ⌜α⌝<n)
```

The same reasoning applies to conjunction

```
maxVar (α ∧ β) with maxVar α | maxVar β
...     | ⌜α⌝ , αpf | ⌜β⌝ , βpf with compare (varidx ⌜α⌝) (varidx ⌜β⌝)
...                              | less ⌜α⌝≤⌜β⌝ = ⌜β⌝ , maxIs⌜β⌝
    where
      maxIs⌜β⌝ : ∀ n → varidx ⌜β⌝ < n → var n FreshIn (α ∧ β)
      maxIs⌜β⌝ n ⌜β⌝<n = αpf n (≤trans (sn≤sm ⌜α⌝≤⌜β⌝) ⌜β⌝<n) ∧ βpf n ⌜β⌝<n
...                              | more ⌜β⌝≤⌜α⌝ = ⌜α⌝ , maxIs⌜α⌝
    where
      maxIs⌜α⌝ : ∀ n → varidx ⌜α⌝ < n → var n FreshIn (α ∧ β)
      maxIs⌜α⌝ n ⌜α⌝<n = αpf n ⌜α⌝<n ∧ βpf n (≤trans (sn≤sm ⌜β⌝≤⌜α⌝) ⌜α⌝<n)
```

and disjunction.

```
maxVar (α ∨ β) with maxVar α | maxVar β
...     | ⌜α⌝ , αpf | ⌜β⌝ , βpf with compare (varidx ⌜α⌝) (varidx ⌜β⌝)
...                              | less ⌜α⌝≤⌜β⌝ = ⌜β⌝ , maxIs⌜β⌝
    where
      maxIs⌜β⌝ : ∀ n → varidx ⌜β⌝ < n → var n FreshIn (α ∨ β)
      maxIs⌜β⌝ n ⌜β⌝<n = αpf n (≤trans (sn≤sm ⌜α⌝≤⌜β⌝) ⌜β⌝<n) ∨ βpf n ⌜β⌝<n
...                              | more ⌜β⌝≤⌜α⌝ = ⌜α⌝ , maxIs⌜α⌝
    where
      maxIs⌜α⌝ : ∀ n → varidx ⌜α⌝ < n → var n FreshIn (α ∨ β)
      maxIs⌜α⌝ n ⌜α⌝<n = αpf n ⌜α⌝<n ∨ βpf n (≤trans (sn≤sm ⌜β⌝≤⌜α⌝) ⌜α⌝<n)
```

For a universal generalisation $\forall x \alpha$, take the greater of $\lceil \alpha \rceil$ and $x$.

```
maxVar (Λ x α) with maxVar α
...                | ⌜α⌝ , αpf with compare (varidx x) (varidx ⌜α⌝)
...                            | less x≤⌜α⌝ = ⌜α⌝ , maxIs⌜α⌝
  where
      maxIs⌜α⌝ : ∀ n → varidx ⌜α⌝ < n → var n FreshIn Λ x α
      maxIs⌜α⌝ n ⌜α⌝<n = Λ (λ { refl → ℕdisorder ⌜α⌝<n x≤⌜α⌝ }) (αpf n ⌜α⌝<n)
...                            | more ⌜α⌝≤x = x , maxIsx
  where
      maxIsx : ∀ n → varidx x < n → var n FreshIn Λ x α
      maxIsx n x<n = Λ (λ { refl → ℕdisorder x<n ≤refl })
                         (αpf n (≤trans (sn≤sm ⌜α⌝≤x) x<n))
```

The same applies for existential generalisation.

```
maxVar (V x α) with maxVar α
...                | ⌜α⌝ , αpf with compare (varidx x) (varidx ⌜α⌝)
...                            | less x≤⌜α⌝ = ⌜α⌝ , maxIs⌜α⌝
```

```
  where
      maxIs⌜α⌝ : ∀ n → varidx ⌜α⌝ < n → var n FreshIn V x α
      maxIs⌜α⌝ n ⌜α⌝<n = V (λ { refl → ℕdisorder ⌜α⌝<n x≤⌜α⌝ }) (apf n ⌜α⌝<n)
...                            | more ⌜α⌝≤x = x , maxIsx
  where
      maxIsx : ∀ n → varidx x < n → var n FreshIn V x α
      maxIsx n x<n = V (λ { refl → ℕdisorder x<n ≤refl })
                         (apf n (≤trans (sn≤sm ⌜α⌝≤x) x<n))
```
□

Finally, a fresh variable can be extracted by choosing the successor of the variable given by the proof above.

```
fresh : ∀ α → Σ Variable (_FreshIn α)
fresh α with maxVar α
...      | ⌜α⌝ , apf = var (suc (varidx ⌜α⌝)) , apf (suc (varidx ⌜α⌝)) ≤refl
```

# 7 Ensemble.lagda

Serious consideration must be given to the data type used to describe the context of a natural deduction tree. In a proof tree for $\Gamma \vdash \alpha$, it must be verified that the remaining open assumptions are all members of $\Gamma$, so the type must have a notion of 'subset'. For universal generalisation introduction, and existential generalisation elimination, it will also be necessary to verify that a given variable is not free in any open assumption, so the type must also have a notion for a predicate holding on all elements. Throughout the natural deduction proof, the collection of open assumptions is modified, either by making new assumptions, by combining collections of assumptions, or by discharging assumptions. Finally, while we will be giving proofs about natural deduction trees, we would also like to give proofs regarding actual formulae (and axiom schemes). Giving natural deduction proofs in this system should correspond closely to doing natural deduction (from the bottom up) by hand. There should not be any need for operations other than the usual rules for natural deduction (with a single exception at the beginning of the proof, as will be shown later). Any manipulation of the context should be done automatically by Agda, and proofs regarding variable freedom and open assumptions should be solvable using Agda's proof search.

The `List` (or `Vec`) type is not suitable. While removal of elements from a list of formulae can be defined with a function, it is unwieldy to give proofs regarding the results of such computations, as they depend on equality-checking of formulae, and so proofs must include both the case where the equality is as expected, and the degenerate case.[4]

An implementation of classical propositional logic in the style of natural deduction was given in [5]. While this does use (something equivalent to) lists, it requires frequent use of extra deduction rules for weakening the context. This would not be suitable for a natural deduction assistant, and it also does not solve the problems given above for first order logic.

Predicates can be used to store collections of values, in the manner of set comprehension. We define the type `Ensemble` as another name for `Pred`. It will be used to refer to predicates which have been created in a manner to follow. This is only for ease of understanding, and is not an actual restriction. Ensembles will resemble finite sets.

```
Ensemble : Set → Set₁
Ensemble A = A → Set
```

Membership is defined by satisfying the predicate.

[4] Examples of this are included in the appendix.

```
infix 4 _∈_ _∉_

_∈_ : {A : Set} → A → Ensemble A → Set
α ∈ αs = αs α

_∉_ : {A : Set} → A → Ensemble A → Set
α ∉ αs = ¬(α ∈ αs)
```

A sensible definition of subset is $A \subset B$ if $\forall x(x \in A \to x \in B)$. However, some ensembles will be defined using negations. If it is absurd for $x$ to be in $A$ (for example, if $A$ is the empty set), then proving that $x \in B$ can be done by either pattern matching to an empty case, or using the lemma `⊥-elim`. However, Agda's proof search will not do pattern matching inside lambda expressions,[5] and it will not find lemmas unless it is hinted to do so. For convenience, we adopt a minimal logic translation by taking the double negative of the right side of the implication, which solves this issue.[6]

```
infix 4 _⊂_
_⊂_ : {A : Set} → Ensemble A → Ensemble A → Set
αs ⊂ βs = ∀ x → x ∈ αs → ¬(x ∉ βs)
```

The empty ensemble and singleton ensembles are defined in the obvious way.

```
∅ : {A : Set} → Ensemble A
∅ = λ _ → ⊥

⟨_⟩ : {A : Set} → A → Ensemble A
⟨ α ⟩ = λ x → x ≡ α
```

It would be reasonable to define union in terms of a disjoint union type, so that a proof of $x \in A \cup B$ would be either a proof of $x \in A$ or of $x \in B$. However, we want Agda's proof search to fill in proofs regarding subsets. For a proof that $A \cup B \subset A \cup B \cup \varnothing$, we would have to do a case analysis on a proof of $x \in A \cup B$. Instead we define $x \in A \cup B$ using functions, so that pattern matching is not necessary in order to make use of such a proof. One definition involving only functions is $x \in A \cup B \coloneqq x \notin A \to x \in B$. We take the double negative of the right side of the implication, for the same reasons as above.

```
infixr 5 _∪_
_∪_ : {A : Set} → Ensemble A → Ensemble A → Ensemble A
(αs ∪ βs) = λ x → x ∉ αs → ¬(x ∉ βs)
```

Instead of defining a set difference, we define notation for removing a single element from an ensemble. Since ensembles will be used only for finite collections, this is not a limitation. A definition using conjunctions is that $x \in A - a$ means $x \in A$ and $x \neq a$. Translating this to functions gives $x \in A - a \coloneqq \neg(x \in A \to x \equiv a)$. Take the contrapositive of the inner implication.

```
infixl 5 _-_
_-_ : {A : Set} → Ensemble A → A → Ensemble A
(αs - α) = λ x → ¬(x ≢ α → x ∉ αs)
```

These definitions allow subset propositions to be proved without case analysis or `⊥-elim` (EFQ), by adopting functional definitions and using double negations. Moreover, the only quantifier used in the definitions is in the definition of `_⊂_`. Since functions are equivalent to implications, we have translated the notion of subset to a proposition of the form $\forall x A$, where $A$ is a formula in the implicational fragment of minimal logic. This is to be expected, since we wanted the proof terms to be simply typed lambda calculus terms, which is precisely equivalent to minimal logic [11].

[5] As of version 2.6.0.

[6] A catalogue of negative translations can be found in [7]. The translation we use is less complete, as we use only enough negations to make the subset predicate minimally provable.

Subset proofs can now be solved by Agda automatically, with good performance. In the case of all natural deduction proofs to follow, Agda solved the subset proof in less than one second (the default time limit for proof search). Moreover, since the implicational fragment of minimal logic is decidable, there are proof search algorithms which will always find a proof if one exists [15].

Of course, ensembles are just predicates, so they can be created in any way that functions can be created. We can define a type to keep track of the creation of a predicate, to ensure it was created using (something equal to) the functions above. Additionally, the type requires that the predicate is over a type with a decidable equality.

```
data Assembled {A : Set} (eq : Decidable≡ A) : Pred A → Set₁ where
  from∅    : Assembled eq ∅
  from⟨_⟩ : (α : A)   → Assembled eq (⟨ α ⟩)
  from_∪_ : ∀{αs βs} → Assembled eq αs → Assembled eq βs
                      → Assembled eq (αs ∪ βs)
  from_-_ : ∀{αs}     → Assembled eq αs → (α : A) → Assembled eq (αs - α)
```

**Proposition 7.0.1.** *Assembled ensembles have decidable membership.*

*Proof.*

```
dec∈ : {A : Set} {eq : Decidable≡ A} {αs : Ensemble A}
       → (x : A) → Assembled eq αs → Dec (x ∈ αs)
```

Nothing is in the empty ensemble.

```
dec∈             x from∅              = no λ x∈∅ → x∈∅
```

Membership of a singleton is defined by an equality, and so its decidability is just the the decidable equality from `Assembled`.

```
dec∈ {_} {eq} x (from⟨ α ⟩)        = eq x α
```

To check membership for a union, simply check first for membership of the left ensemble, then the right. The lambda expression proofs given here are non-trivial, and difficult to interpret, but can be provided by Agda's proof search.

```
dec∈             x (from Aαs ∪ Aβs) with dec∈ x Aαs
...     | yes x∈αs = yes λ x∉αs _ → x∉αs x∈αs
...     | no  x∉αs with dec∈ x Aβs
...                 | yes x∈βs = yes λ _ x∉βs → x∉βs x∈βs
...                 | no  x∉βs = no  λ x∉αs∪βs → x∉αs∪βs x∉αs x∉βs
```

Finally, in the case of an element being removed, use the decidable equality from `Assembled` to check if the given element was removed, and otherwise check if the given element is in the inner ensemble.

```
dec∈ {_} {eq} x (from Aαs - α)   with eq x α
...     | yes refl = no λ α∈αs-α → α∈αs-α λ α≢α _ → α≢α refl
...     | no x≢α   with dec∈ x Aαs
...                 | yes x∈αs = yes λ x≢α→x∉αs → x≢α→x∉αs x≢α x∈αs
...                 | no  x∉αs = no  λ x∈αs-α
                                    → x∈αs-α (λ _ _
                                          → x∈αs-α (λ _ _
                                                → x∈αs-α (λ _
                                                      → x∉αs)))
```

□

Given an ensemble $A$, a sensible definition for a predicate $P$ holding on every element of $A$ would be $\forall x(x \in A \to Px)$. However, for inductively defined predicates (like `_notFreeIn α` for some $\alpha$), this is not easy to work with, either by hand or using proof search. For example, to prove that the variable $y$ is not free in all members of $\{\forall y Q y\} \cup \{\bot\}$, it would be necessary to show that every member is equal to either $\forall y Q y$ or $\bot$, and only then supply the required constructors for each case. Once again, this requires pattern matching.

Instead, for an assembled ensemble, we give a definition for `All` which utilises the structure of the ensemble, and describes what computation must be performed to check that a predicate holds on all members. To do so, maintain a list of all elements which have been removed from the ensemble.

```
infixr 5 _all∪_

data All_[_\_] {A : Set} (P : Pred A) : Ensemble A → List A → Set₁ where
  all∅     : ∀{xs}        → All P [ ∅ ∖ xs ]
```

$P$ holds on all of a singleton if it holds on the element of the singleton, or else if that element has already been removed.

```
  all⟨_⟩  : ∀{α xs}      → P α           → All P [ ⟨ α ⟩ ∖ xs ]
  all⟨-_⟩ : ∀{α xs}      → α List.∈ xs → All P [ ⟨ α ⟩ ∖ xs ]
```

In the case of a union, $P$ must hold on both sides of the union.

```
  _all∪_  : ∀{αs βs xs} → All P [ αs ∖ xs ] → All P [ βs ∖ xs ]
                          → All P [ αs ∪ βs ∖ xs ]
```

Finally, when an ensemble has been created by removing an element from another, check that $P$ holds on the other ensemble for all values other than the removed one.

```
  all-_   : ∀{αs x xs}   → All P [ αs ∖ x :: xs ] → All P [ αs - x ∖ xs ]
```

Now, $P$ holds on all of $\alpha s$ if it holds according to the above procedure, with the removed element list starting empty.

```
All : {A : Set} → Pred A → Ensemble A → Set₁
All P αs = All P [ αs ∖ [] ]
```

**Proposition 7.0.2.** *The definition of* `All` *for assembled ensembles is weaker than the usual set definition.*

*Proof.* We use a lemma to show that this is the case for all values of the removed list of elements, and apply it to the case of the empty list.

```
fAll→All : {A : Set} {eq : Decidable≡ A} {P : Pred A} {αs : Ensemble A}
            → Assembled eq αs → (∀ x → x ∈ αs → P x) → All P αs
fAll→All {A} {eq} {P} Aαs fall = ψ Aαs [] (λ x x∈αs _ → fall x x∈αs)
  where
    ψ : ∀{αs} → Assembled eq αs → ∀ xs
                → (∀ x → x ∈ αs → x List.∉ xs → P x) → All P [ αs ∖ xs ]
    ψ from∅             xs fall∅      = all∅
```

For a singleton $\{\alpha\}$, either $\alpha$ has been removed, or otherwise it has not been removed, in which case we use the functional all to prove $P\alpha$.

```
    ψ from⟨ α ⟩         xs fall⟨α⟩   with List.dec∈ eq α xs
    ...                                 | yes α∈xs = all⟨- α∈xs ⟩
    ...                                 | no  α∉xs = all⟨ fall⟨α⟩ α refl α∉xs ⟩
```

Since unions are defined using a double negation, to show that the functional all for a union means functional all for each of the two ensembles, use a contradiction for each.

```
ψ (from Aαs ∪ Aβs) xs fallαs∪βs = (ψ Aαs xs fallαs)
                                   all∪ (ψ Aβs xs fallβs)
  where
    fallαs : _
    fallαs x x∈αs = fallαs∪βs x (λ x∉αs _     → x∉αs x∈αs)
    fallβs : _
    fallβs x x∈βs = fallαs∪βs x (λ _     x∉βs → x∉βs x∈βs)
```

In the case of $\alpha s - \alpha$, we show first that if $x \in \alpha s$ then $x \in \alpha s - \alpha$, and that if $x \notin \alpha :: xs$ then $x \notin xs$.

```
ψ (from Aαs - α)   xs fallαs-α  = all- (ψ Aαs (α :: xs) fallαs)
  where
    fallαs : _
    fallαs x x∈αs x∉α::xs =
      let x∈αs-α : _
          x∈αs-α x≢α→x∉αs = x≢α→x∉αs (λ x≡α → x∉α::xs List.[ x≡α ]) x∈αs
          x∉xs   : x List.∉ xs
          x∉xs   x∈xs = x∉α::xs (α :: x∈xs)
      in  fallαs-α x x∈αs-α x∉xs                                         □
```

The converse cannot be proved; `All` is in fact strictly weaker than the functional definition. While it could be expected that pattern matching on both `All` and `Assembled` would lead to a proof, this will not work because Agda cannot unify function types. For example, in the case that an ensemble was assembled by `from Aαs ∪ Aβs`, case analysis of the proof of `All P (αs ∪ βs)` does not show that the only constructor is `_all∪_`; Agda cannot determine that `λ x → x ∉ αs → ¬(x ∉ βs)` does not unify with `λ _ → ⊥`, so `all∅` may or may not be a constructor. If we wanted a stronger type which is equivalent to the functional definition, the assembled structure would need to be included in `All`.

We can use the `All` predicate to define the restriction that certain deductions are valid only if a given variable is not free in an ensemble of open assumptions. For the usual use case (i.e. cases other than abstract proof tree manipulation where variable freedom is determined by some lemma), Agda's proof search will find the required proof. However, due to the above limitations with unification of functions, Agda does not see that there is only one constructor for each non-singleton ensemble, so the search algorithm is not fast. For larger proof trees, it is necessary to increase the timeout from the default one second to ten seconds. This could also be resolved by including the assembled structure in `All`.

## 8 Deduction.lagda

We now define the type of natural deductions, using the deduction rules of of [10]. Given $\Gamma$ and $\alpha$, anything that the type checker confirms as being of type $\Gamma \vdash \alpha$ is a valid natural deduction proof of $\alpha$ from assumptions $\Gamma$, and so is a proof of $\alpha$ from $\Gamma$ over minimal logic.

First, some shorthand.

```
private
  _NotFreeInAll_ : Variable → Ensemble Formula → Set₁
  x NotFreeInAll Γ = All (x NotFreeIn_) Γ
```

Now for the natural deduction rules.

```
infix 1 _⊢_ ⊢_

data _⊢_ : Ensemble Formula → Formula → Set₁ where
```

The first constructor is not a deduction rule, in that it does not change the type of the deduction. It will be used for typesetting later, for abbreviating a previously proved deduction from no assumptions. This will be used for lemmas, and for applying assumed axiom schemes.

```
  cite         : ∀{α} → String → ∅ ⊢ α → ∅ ⊢ α
```

The following constructor exists primarily to 'normalise' Γ, for example replacing a proof of $\{\alpha\} - \alpha \vdash \beta$ with a proof of $\varnothing \vdash \beta$. It is also necessary for weakening results, for example from $\Gamma \vdash \alpha$ to $\Gamma, \beta \vdash \alpha$. While this is not one of the usual deduction rules, it will need to be used only at the beginning of a proof to finalise the ensemble of assumptions. We require that an assembled ensemble is given, so that membership remains decidable.

```
  close        : ∀{Γ Δ α} → Assembled formulaEq Δ → Γ ⊂ Δ → Γ ⊢ α → Δ ⊢ α
```

The remaining constructors correspond precisely to the usual natural deduction rules. Agda's comment syntax (--) allows these rules to be formatted as Gentzen-style inferences.

```
  assume       : (α : Formula)
                 →                                ⟨ α ⟩ ⊢ α

  arrowintro   : ∀{Γ β} → (α : Formula)
                 →                                  Γ ⊢ β
                                                ---------------- ⇒⁺
                 →                               Γ - α ⊢ α ⇒ β

  arrowelim    : ∀{Γ₁ Γ₂ α β}
                 →                    Γ₁ ⊢ α ⇒ β     →     Γ₂ ⊢ α
                                    ------------------------------ ⇒⁻
                 →                           Γ₁ ∪ Γ₂ ⊢ β

  conjintro    : ∀{Γ₁ Γ₂ α β}
                 →                      Γ₁ ⊢ α     →     Γ₂ ⊢ β
                                      -------------------------- ∧⁺
                 →                        Γ₁ ∪ Γ₂ ⊢ α ∧ β

  conjelim     : ∀{Γ₁ Γ₂ α β γ}
                 →                    Γ₁ ⊢ α ∧ β     →     Γ₂ ⊢ γ
                                    ------------------------------ ∧⁻
                 →                      Γ₁ ∪ (Γ₂ - α - β) ⊢ γ

  disjintro₁   : ∀{Γ α} → (β : Formula)
                 →                                  Γ ⊢ α
                                                  ------------ ∨⁺₁
                 →                                 Γ ⊢ α ∨ β

  disjintro₂   : ∀{Γ β} → (α : Formula)
                 →                                  Γ ⊢ β
                                                  ------------ ∨⁺₂
                 →                                 Γ ⊢ α ∨ β

  disjelim     : ∀{Γ₁ Γ₂ Γ₃ α β γ}
                 →             Γ₁ ⊢ α ∨ β     →     Γ₂ ⊢ γ     →     Γ₃ ⊢ γ
                             -------------------------------------------- ∨⁻
                 →                    Γ₁ ∪ (Γ₂ - α) ∪ (Γ₃ - β) ⊢ γ
```

The constructors for first order logic require an extra proof to be supplied, either of variable freedom or variable substitution. The propositions proved here have been formulated so that Agda's built-in proof search should be able to supply them.

```
  univintro    : ∀{Γ α} → (x : Variable)
                 → x NotFreeInAll Γ
                 →                              Γ ⊢ α
                                              ----------- ∀⁺
                 →                            Γ ⊢ Λ x α

  univelim     : ∀{Γ α x α[x/t]} → (t : Term)
                 → α [ x / t ]≡ α[x/t]
                 →                            Γ ⊢ Λ x α
                                              ------------ ∀⁻
                 →                            Γ ⊢ α[x/t]

  existintro   : ∀{Γ α α[x/t]} → (t : Term) → (x : Variable)
                 → α [ x / t ]≡ α[x/t]
                 →                            Γ ⊢ α[x/t]
                                              ------------ ∃⁺
                 →                            Γ ⊢ V x α

  existelim    : ∀{Γ₁ Γ₂ α β x}
                 → x NotFreeInAll (⟨ β ⟩ ∪ (Γ₂ - α))
                 →                      Γ₁ ⊢ V x α    →    Γ₂ ⊢ β
                                      --------------------------- ∃⁻
                 →                          Γ₁ ∪ (Γ₂ - α) ⊢ β
```

Finally, we define the following shorthand.

```
⊢_ : Formula → Set₁
⊢ α = ∅ ⊢ α
```

It is trivial to show that the context of a deduction is assembled (and so membership is decidable), simply by recursing over the deduction rules. The proof is omitted.

```
assembled-context : ∀{Γ α} → Γ ⊢ α → Assembled formulaEq Γ
-- Proof omitted.
```

# 9 Formula equivalence

Avoiding conflicts of variables with the same name can be done algorithmically by machines (including internally by Agda) by using a nameless notation. Since natural deduction is intended to be used by humans, we prefer to state that formulae can be used equivalently if they are equivalent up to renaming of bound variables [4, 10]. This has not been included in the definition of formulae or of natural deduction. To do so would introduce an extra complication to the deduction rules, as every step in a natural deduction proof would have to include a proof that the conclusion is equivalent to the desired one. However, it is possible to prove that this is unnecessary; in the system given, if $\Gamma \vdash \alpha$, and $\alpha$ is equivalent to $\alpha'$ up to the renaming of bound variables, then $\Gamma \vdash \alpha'$.

## 9.1 Formula equivalence

Formulae are *equivalent* if they are equal up to renaming bound variables.[7] Here, renaming means substituting a variable for another variable which is not free, so that the meaning of the formula does not change.

```
infix 50 _≈_
data _≈_ : Formula → Formula → Set where
```

First, the trivial cases for equivalence, coming from equivalence of components.

```
  atom : ∀ r ts → atom r ts ≈ atom r ts
  _⇒_  : ∀{α β α′ β′} → α ≈ α′ → β ≈ β′ → α ⇒ β ≈ α′ ⇒ β′
  _∧_  : ∀{α β α′ β′} → α ≈ α′ → β ≈ β′ → α ∧ β ≈ α′ ∧ β′
  _∨_  : ∀{α β α′ β′} → α ≈ α′ → β ≈ β′ → α ∨ β ≈ α′ ∨ β′
  Λ    : ∀{α α′} → ∀ x → α ≈ α′ → Λ x α ≈ Λ x α′
  V    : ∀{α α′} → ∀ x → α ≈ α′ → V x α ≈ V x α′
```

Now, the case for renaming the quantifying variable of a generalisation. The resulting component must also be replaceable with an equivalent component, as other bound variable renaming may occur inside.

```
  Λ/   : ∀{α β β′ x y} → y NotFreeIn α → α [ x / varterm y ]≡ β
                       → β ≈ β′ → Λ x α ≈ Λ y β′
  V/   : ∀{α β β′ x y} → y NotFreeIn α → α [ x / varterm y ]≡ β
                       → β ≈ β′ → V x α ≈ V y β′
```

For equivalence to be symmetric, the following dual form of bound variable renaming must be derivable.

```
  Λ/′  : ∀{α α′ β′ x y} → α ≈ α′ → y NotFreeIn α′
                        → α′ [ x / varterm y ]≡ β′ → Λ x α ≈ Λ y β′
  V/′  : ∀{α α′ β′ x y} → α ≈ α′ → y NotFreeIn α′
                        → α′ [ x / varterm y ]≡ β′ → V x α ≈ V y β′
```

It may be that the latter two rules are derivable. However, if this is so, proving this would require several lemmas which will be otherwise unnecessary. As the goal here is to prove that equivalent formulae are equivalently derivable, having extra constructors for equivalence will not weaken this result. It will be shown later that it would be more ‘natural’ to adopt the rules `Λ/` and `V/′` and to derive `Λ/′` and `V/` if possible.

**Lemma 9.1.1.** *Formula equivalence is symmetric.*

*Proof.*

```
≈sym : ∀{α α′} → α ≈ α′ → α′ ≈ α
```

For the trivial definitions, the proof is similarly trivial.

```
≈sym (atom r ts)   = atom r ts
≈sym (α≈α′ ⇒ β≈β′) = ≈sym α≈α′ ⇒ ≈sym β≈β′
≈sym (α≈α′ ∧ β≈β′) = ≈sym α≈α′ ∧ ≈sym β≈β′
≈sym (α≈α′ ∨ β≈β′) = ≈sym α≈α′ ∨ ≈sym β≈β′
≈sym (Λ x α≈α′)    = Λ x (≈sym α≈α′)
≈sym (V x α≈α′)    = V x (≈sym α≈α′)
```

In the case of bound variable renaming, the dual constructor is used. If the bound variable is being renamed with itself, then the previous trivial proof is given instead.

[7] This is often called $\alpha$-equivalence.

```
≈sym {Λ x α} {Λ y β'} (Λ/ y∉α α[x/y]≡β β≈β') with varEq x y
...     | yes refl rewrite subIdentFunc α[x/y]≡β = Λ x (≈sym β≈β')
...     | no  x≢y  = Λ/' (≈sym β≈β') (subNotFree (varterm x≢y) α[x/y]≡β)
                        (subInverse y∉α α[x/y]≡β)
≈sym {V x α} {V y β'} (V/ y∉α α[x/y]≡β β≈β') with varEq x y
...     | yes refl rewrite subIdentFunc α[x/y]≡β = V x (≈sym β≈β')
...     | no  x≢y  = V/' (≈sym β≈β') (subNotFree (varterm x≢y) α[x/y]≡β)
                        (subInverse y∉α α[x/y]≡β)
≈sym {Λ x α} {Λ y β'} (Λ/' α≈α' y∉α' α'[x/y]≡β') with varEq x y
...     | yes refl rewrite subIdentFunc α'[x/y]≡β' = Λ x (≈sym α≈α')
...     | no  x≢y  = Λ/ (subNotFree (varterm x≢y) α'[x/y]≡β')
                       (subInverse y∉α' α'[x/y]≡β') (≈sym α≈α')
≈sym {V x α} {V y β'} (V/' α≈α' y∉α' α'[x/y]≡β') with varEq x y
...     | yes refl rewrite subIdentFunc α'[x/y]≡β' = V x (≈sym α≈α')
...     | no  x≢y  = V/ (subNotFree (varterm x≢y) α'[x/y]≡β')
                       (subInverse y∉α' α'[x/y]≡β') (≈sym α≈α')
```
□

If a variable is not free in $\alpha$, then it should not be free in $\alpha[x/t]$, assuming that the variable does not appear in $t$. The proof simply comes from the definition of variable substitution and freedom for terms.

```
notFreeSub : ∀{α β x t z} → z NotFreeIn α → z NotInTerm t
             → α [ x / t ]≡ β → z NotFreeIn β
-- Proof omitted.
```

Variable freedom is preserved by formula equivalence. This is proved using the lemma above, noting that if $z$ is bound in $\alpha$, then it is either also bound in $\alpha'$ or else has been renamed and so does not appear, and if $z$ does not appear in $\alpha$ then it either also does not appear in $\alpha'$ or some bound variable has been renamed to it, so it is bound.

```
≈notFree : ∀{α α' z} → α ≈ α' → z NotFreeIn α → z NotFreeIn α'
-- Proof omitted.
```

## 9.2 Deriving the rename rule

We want to derive the deduction rule $\alpha \approx \alpha' \to \Gamma \vdash \alpha \to \Gamma \vdash \alpha'$. As $\Gamma$ is the same in both deductions, changing variable names within the deduction of $\Gamma \vdash \alpha$ will not suffice. Instead, the deduction tree is extended to obtain the new conclusion.

Proving this rule has a termination issue for the case that $\alpha$ is an implication, which will be explained below. However, it is possible to prove the stronger rule $\alpha \approx \alpha' \to (\Gamma \vdash \alpha \leftrightarrow \Gamma \vdash \alpha')$. A simplified notion of '$\leftrightarrow$' in Agda can be defined as follows.[8]

```
private
  record _↔_ (A : Set₁) (B : Set₁) : Set₁ where
    field
      ⟨→⟩ : A → B
      ⟨←⟩ : B → A
  open _↔_
```

We can now define the stronger rename rule.

```
renameIff : ∀{Γ α α'} → α ≈ α' → (Γ ⊢ α) ↔ (Γ ⊢ α')
```

Clearly, the rename rule can be derived from `renameIff`.

[8] This is simplified in that it requires all types involved to be of type `Set₁`, which is enough for our purposes.

```
rename       : ∀{Γ α α′}
               → α ≈ α′
               →                                 Γ ⊢ α
                                                ---------
               →                                 Γ ⊢ α′
rename α≈α′ = ⟨→⟩ (renameIff α≈α′)
```

It remains to prove `renameIff`. In the natural deduction proofs that follow, the subset proofs for `close` were found automatically by Agda's proof search. The variable freedom proofs were also found, except where extra reasoning (regarding substitution and equivalence lemmas) were required. Such lemmas are necessary when manipulating proof trees in the abstract.

*Proof.* The atomic case is trivial, since an atomic formula is equivalent only to itself.

```
⟨→⟩ (renameIff {Γ} {atom r ts} {.(atom r ts)} (atom .r .ts)) d = d
⟨←⟩ (renameIff {Γ} {atom r ts} {.(atom r ts)} (atom .r .ts)) d = d
```

The proof tree for the implication case is extended as follows.

```
⟨→⟩ (renameIff {Γ} {α ⇒ β} {α′ ⇒ β′} (α≈α′ ⇒ β≈β′)) Γ⊢α⇒β =
```

$$\cfrac{\cfrac{\cfrac{\begin{matrix}\Gamma\\ \vdots\\ \alpha\to\beta\end{matrix} \qquad \cfrac{[\alpha']}{\alpha}\,\text{induction}}{\beta}}{\beta'}\,\text{induction}}{\alpha'\to\beta'}\,\to^-$$

One of the induction steps involves invoking the rename rule on $\alpha' \approx \alpha$ and the assumption of $\alpha'$. We have $\alpha \approx \alpha'$, and `≈sym` shows that formula equivalence is symmetric. However, calling `rename` on `≈sym α≈α′` would not be structurally recursive, because Agda cannot determine that `≈sym α≈α′` is structurally smaller than `α≈α′ ⇒ β≈β′`. This is the reason for proving `renameIff` instead of proving `rename` directly; we have access to a proof of $\alpha' \vdash \alpha$ by using the opposite direction of `renameIff`.

```
  close
   (assembled-context Γ⊢α⇒β)
   (λ x z z₁ → z (λ z₂ z₃ → z₃ z₁ z₂))
   (arrowintro α′
    (⟨→⟩ (renameIff β≈β′)
     (arrowelim
      Γ⊢α⇒β
      -- `rename (≈sym α≈α′) (assume α′)` would not be structurally recursive
      (⟨←⟩ (renameIff α≈α′)
       (assume α′)))))
```

The other direction has the same proof, with $\alpha$ swapped with $\alpha'$, $\beta$ swapped with $\beta'$, and the opposite directions of `renameIff` used.

```
⟨←⟩ (renameIff {Γ} {α ⇒ β} {α′ ⇒ β′} (α≈α′ ⇒ β≈β′)) Γ⊢α′⇒β′ =
-- Proof omitted.
```

The proof tree for the conjunction case is extended as follows.

```
⟨→⟩ (renameIff {Γ} {α ∧ β} {α′ ∧ β′} (α≈α′ ∧ β≈β′)) Γ⊢α∧β =
```

$$\dfrac{\begin{array}{c}\Gamma\\ \vdots\\ \alpha\wedge\beta\end{array} \qquad \dfrac{\dfrac{[\alpha]}{\alpha'}\,\text{induction} \qquad \dfrac{[\beta]}{\beta'}\,\text{induction}}{\alpha'\wedge\beta'}\,\wedge^+}{\alpha'\wedge\beta'}\,\wedge^-$$

```
  close
   (assembled-context Γ⊢α∧β)
   (λ x z z₁ → z z₁ (λ z₂ → z₂ (λ z₃ z₄ → z₄ (λ z₅ z₆ → z₆ z₅ z₃))))
   (conjelim
    Γ⊢α∧β
    (conjintro
     ((⟨→⟩ (renameIff α≈α')
      (assume α))
     ((⟨→⟩ (renameIff β≈β')
      (assume β))))
```

Again, the other direction is obtained by reversing the use of equivalences.

```
⟨←⟩ (renameIff {Γ} {α ∧ β} {α' ∧ β'} (α≈α' ∧ β≈β')) Γ⊢α'∧β' =
-- Proof omitted.
```

The proof tree for the disjunction case is extended as follows.

```
⟨→⟩ (renameIff {Γ} {α ∨ β} {α' ∨ β'} (α≈α' ∨ β≈β')) Γ⊢α∨β =
```

$$\dfrac{\begin{array}{c}\Gamma\\ \vdots\\ \alpha\vee\beta\end{array} \qquad \dfrac{\dfrac{[\alpha]}{\alpha'}\,\text{induction}}{\alpha'\vee\beta'}\,\vee^+ \qquad \dfrac{\dfrac{[\beta]}{\beta'}\,\text{induction}}{\alpha'\vee\beta'}\,\vee^+}{\alpha'\vee\beta'}\,\vee^-$$

```
  close
   (assembled-context Γ⊢α∨β)
   (λ x z z₁
    → z z₁ (λ z₂ → z₂ (λ z₃ → z₃ (λ z₄ → z₄)) (λ z₃ → z₃ (λ z₄ → z₄))))
   (disjelim
    Γ⊢α∨β
    (disjintro₁ β'
     ((⟨→⟩ (renameIff α≈α')
      (assume α)))
    (disjintro₂ α'
     ((⟨→⟩ (renameIff β≈β')
      (assume β))))
```

Again, the other direction is obtained by reversing the use of equivalences.

```
⟨←⟩ (renameIff {Γ} {α ∨ β} {α' ∨ β'} (α≈α' ∨ β≈β')) Γ⊢α'∨β' =
-- Proof omitted.
```

The first case for universal generalisation is where the bound variable is not renamed.

```
⟨→⟩ (renameIff {Γ} {Λ x α} {Λ .x α'} (Λ y α≈α')) Γ⊢∀xα =
```

Since $x$ may be free $\Gamma$, we use arrow introduction and elimination so that $\Gamma$ is not assumed when the universal generalisation is re-introduced.

$$\dfrac{\dfrac{\dfrac{\dfrac{\dfrac{[\forall x\alpha]}{\alpha}\,\forall^-}{\alpha'}\,\text{induction}}{\forall x\alpha'}\,\forall^+}{\forall x\alpha \to \forall x\alpha'}\,\to^+ \qquad \begin{matrix}\Gamma\\ \vdots\\ \forall x\alpha\end{matrix}}{\forall x\alpha'}\,\to^-$$

```
  close
   (assembled-context Γ⊢∀xα)
   (λ x z → z (λ z₁ → z₁ (λ z₂ → z₂)))
   (arrowelim
    (arrowintro (Λ x α)
     (univintro x (all⟨ Λ↓ x α ⟩)
      (⟨→⟩ (renameIff α≈α′)
       (univelim (varterm x) (ident α x)
        (assume (Λ x α))))))
    Γ⊢∀xα)
```

Again, the other direction is obtained by reversing the use of equivalences.

```
⟨←⟩ (renameIff {Γ} {Λ x α} {Λ .x α′} (Λ y α≈α′)) Γ⊢∀xα′ =
-- Proof omitted.
```

The second case for universal generalisation renames the bound variable, then follows another equivalence.

```
⟨→⟩ (renameIff {Γ} {Λ x α} {Λ y β′} (Λ/ y∉α α[x/y]≡β β≈β′)) Γ⊢∀xα =
```

Since $x$ is being renamed to $y$, we know $y$ is not free in $\alpha$, and so it is also not free in $\forall x\alpha$. Define $\beta \coloneqq \alpha[x/y]$.

$$\dfrac{\dfrac{\dfrac{\dfrac{\dfrac{[\forall x\alpha]}{\beta}\,\forall^-}{\beta'}\,\text{induction}}{\forall y\beta'}\,\forall^+}{\forall x\alpha \to \forall y\beta'}\,\to^+ \qquad \begin{matrix}\Gamma\\ \vdots\\ \forall x\alpha\end{matrix}}{\forall y\beta'}\,\to^-$$

```
  close
   (assembled-context Γ⊢∀xα)
   (λ x₁ z → z (λ z₁ → z₁ (λ z₂ → z₂)))
   (arrowelim
    (arrowintro (Λ x α)
     (univintro y all⟨ Λ x y∉α ⟩
      (⟨→⟩ (renameIff β≈β′)
       (univelim (varterm y) α[x/y]≡β
        (assume (Λ x α))))))
    Γ⊢∀xα)
```

The other direction varies depending on if $x$ is equal to $y$.

```
⟨←⟩ (renameIff {Γ} {Λ x α} {Λ y β′} (Λ/ y∉α α[x/y]≡β β≈β′)) Γ⊢∀yβ′
    with varEq x y
```

In the degenerate case where $x = y$, we have $\beta = \alpha$.

$$\dfrac{\dfrac{\dfrac{\dfrac{\dfrac{[\forall x\beta']}{\beta'}\,\forall^-}{\alpha}\,\text{induction}}{\forall x\alpha}\,\forall^+}{\forall x\beta' \to \forall x\alpha}\,\to^+ \qquad \begin{matrix}\Gamma\\ \vdots\\ \forall x\beta'\end{matrix}}{\forall x\alpha}\,\to^-$$

```
... | yes refl rewrite subIdentFunc α[x/y]≡β =
  close
   (assembled-context Γ⊢∀yβ')
   (λ x z → z (λ z₁ → z₁ (λ z₂ → z₂)))
   (arrowelim
    (arrowintro (Λ x β')
     (univintro x all⟨ Λ↓ x β' ⟩
      (⟨←⟩ (renameIff β≈β')
       (univelim (varterm x) (ident β' x)
        (assume (Λ x β'))))))
    Γ⊢∀yβ')
```

Otherwise, $\beta[y/x] = \alpha$, and $x$ is not free in $\forall y\beta'$ because $x$ is not free in $\beta$, since $\beta$ is obtained by substituting $x$ with $y$ in $\alpha$.

$$\dfrac{\dfrac{\dfrac{\dfrac{\dfrac{\dfrac{\dfrac{[\forall y\beta']}{\beta'}\,\forall^-}{\beta}\,\text{induction}}{\forall y\beta}\,\forall^+}{\alpha}\,\forall^-}{\forall x\alpha}\,\forall^+}{\forall y\beta' \to \forall x\alpha}\,\to^+ \qquad \begin{matrix}\Gamma\\ \vdots\\ \forall y\beta'\end{matrix}}{\forall x\alpha}\,\to^-$$

```
... | no  x≢y  =
  close
   (assembled-context Γ⊢∀yβ')
   (λ x z → z (λ z₁ → z₁ (λ z₂ → z₂)))
   (arrowelim
    (arrowintro (Λ y β')
     (univintro x
      all⟨ Λ y (≈notFree β≈β' (subNotFree (varterm x≢y) α[x/y]≡β)) ⟩
      (univelim (varterm x) (subInverse y∉α α[x/y]≡β)
       (univintro y all⟨ Λ↓ y β' ⟩
        (⟨←⟩ (renameIff β≈β')
         (univelim (varterm y) (ident β' y)
          (assume (Λ y β'))))))))
    Γ⊢∀yβ')
```

The third case is the dual of the second.

```
⟨→⟩ (renameIff {Γ} {Λ x α} {Λ y β'} (Λ/' α≈α' y∉α' α'[x/y]≡β')) Γ⊢∀xα =
```

$$\dfrac{\dfrac{\dfrac{\dfrac{\dfrac{\dfrac{\dfrac{[\forall x\alpha]}{\alpha}\forall^-}{\alpha'}\text{induction}}{\forall x\alpha'}\forall^+}{\beta'}\forall^-}{\forall y\beta'}\forall^+}{\forall x\alpha \to \forall y\beta'}\to^+ \qquad \begin{matrix}\Gamma\\ \vdots\\ \forall x\alpha\end{matrix}}{\forall y\beta'}\to^-$$

```
  close
   (assembled-context Γ⊢∀xα)
   (λ x₁ z → z (λ z₁ → z₁ (λ z₂ → z₂)))
   (arrowelim
    (arrowintro (Λ x α)
     (univintro y all⟨ ≈notFree (Λ x (≈sym α≈α')) (Λ x y∉α') ⟩
      (univelim (varterm y) α'[x/y]≡β'
       (univintro x all⟨ Λ↓ x α ⟩
        (⟨→⟩ (renameIff α≈α')
         (univelim (varterm x) (ident α x)
          (assume (Λ x α))))))))
    Γ⊢∀xα)
```

The other direction varies depending on if $x$ is equal to $y$.

```
⟨←⟩ (renameIff {Γ} {Λ x α} {Λ y β'} (Λ/' α≈α' y∉α' α'[x/y]≡β')) Γ⊢∀yβ'
    with varEq x y
```

In the degenerate case where $x = y$, we have $\alpha' = \beta'$.

$$\dfrac{\dfrac{\dfrac{\dfrac{\dfrac{[\forall x\beta']}{\beta'}\forall^-}{\alpha}\text{induction}}{\forall x\alpha}\forall^+}{\forall x\beta' \to \forall x\alpha}\to^+ \qquad \begin{matrix}\Gamma\\ \vdots\\ \forall x\beta'\end{matrix}}{\forall x\alpha}\to^-$$

```
... | yes refl rewrite subIdentFunc α'[x/y]≡β' =
  close
   (assembled-context Γ⊢∀yβ')
   (λ x z → z (λ z₁ → z₁ (λ z₂ → z₂)))
   (arrowelim
    (arrowintro (Λ x β')
     (univintro x all⟨ Λ↓ x β' ⟩
      (⟨←⟩ (renameIff α≈α')
       (univelim (varterm x) (ident β' x)
        (assume (Λ x β'))))))
    Γ⊢∀yβ')
```

Otherwise, $\beta'[y/x] = \alpha'$, and $x$ is not free in $\forall y\beta'$ since $\beta'$ has been obtained by substituting $x$ with $y$ in $\alpha'$.

$$\dfrac{\dfrac{\dfrac{\dfrac{\dfrac{[\forall y\beta']}{\alpha'}\forall^-}{\alpha}\text{induction}}{\forall x\alpha}\forall^+}{\forall y\beta' \to \forall x\alpha}\to^+ \qquad \begin{matrix}\Gamma\\ \vdots\\ \forall y\beta'\end{matrix}}{\forall x\alpha}\to^-$$

```
... | no  x≢y  =
  close
   (assembled-context Γ⊢∀yβ′)
   (λ x z → z (λ z₁ → z₁ (λ z₂ → z₂)))
   (arrowelim
    (arrowintro (Λ y β′)
     (univintro x all⟨ Λ y (subNotFree (varterm x≢y) α′[x/y]≡β′) ⟩
      ((⟨←⟩ (renameIff α≈α′)
       (univelim (varterm x) (subInverse y∉α′ α′[x/y]≡β′)
        (assume (Λ y β′)))))))
    Γ⊢∀yβ′)
```

Finally, we examine the existential generalisation cases. The first case is where the bound variable is not renamed.

```
⟨→⟩ (renameIff {Γ} {V x α} {V .x α′} (V y α≈α′)) Γ⊢∃xα =
```

$$\dfrac{\begin{matrix}\Gamma\\ \vdots\\ \exists x\alpha\end{matrix} \qquad \dfrac{\dfrac{[\alpha]}{\alpha'}\,\text{induction}}{\exists x\alpha'}\,\exists^+}{\exists x\alpha'}\,\exists^-$$

```
  close
   (assembled-context Γ⊢∃xα)
   (λ x z z₁ → z z₁ (λ z₂ → z₂ (λ z₃ → z₃)))
   (existelim (all⟨ V↓ x α′ ⟩ all∪ (all- all⟨- [ refl ] ⟩))
    Γ⊢∃xα
    (existintro (varterm x) x (ident α′ x)
     ((⟨→⟩ (renameIff α≈α′)
      (assume α))))
```

The reverse direction is the same, with equivalences reversed.

```
⟨←⟩ (renameIff {Γ} {V x α} {V .x α′} (V y α≈α′)) Γ⊢∃xα′ =
-- Proof omitted.
```

The second case for existential generalisation renames the bound variable, then follows another equivalence. The proof depends on whether $x$ is equal to $y$.

```
⟨→⟩ (renameIff {Γ} {V x α} {V y β′} (V/ y∉α α[x/y]≡β β≈β′)) Γ⊢∃xα
    with varEq x y
```

Since $\beta = \alpha[x/y]$, we have $\alpha = \beta[y/x]$. If $x \neq y$, then $x$ cannot be free in $\beta$, and so it is also not free in $\exists y\beta$.

$$\dfrac{\begin{matrix}\Gamma\\ \vdots\\ \exists x\alpha\end{matrix} \qquad \dfrac{\dfrac{[\alpha]}{\exists y\beta}\,\exists^+ \qquad \dfrac{\dfrac{[\beta]}{\beta'}\,\text{induction}}{\exists y\beta'}\,\exists^+}{\exists y\beta'}\,\exists^-}{\exists y\beta'}\,\exists^-$$

```
... | no  x≢y  =
  close
   (assembled-context Γ⊢∃xα)
   (λ x₁ z z₁ → z z₁ (λ z₂ → z₂ (λ z₃ z₄ → z₄ z₃ (λ z₅ → z₅ (λ z₆ → z₆)))))
```

```
  (existelim (all⟨ V y (≈notFree β≈β' (subNotFree (varterm x≢y) α[x/y]≡β)) ⟩
              all∪ (all- (all⟨- [ refl ] ⟩ all∪ (all- all⟨- [ refl ] ⟩))))
   Γ⊢∃xα
   (existelim (all⟨ V↓ y β' ⟩ all∪ (all- all⟨- [ refl ] ⟩))
    (existintro (varterm x) y (subInverse y∉α α[x/y]≡β)
     (assume α))
    (existintro (varterm y) y (ident β' y)
     (⟨→⟩ (renameIff β≈β')
      (assume _)))))
```

In the degenerate case, we have $\beta = \alpha$.

$$\dfrac{\begin{matrix}\Gamma \\ \vdots \\ \exists x\alpha\end{matrix} \qquad \dfrac{\dfrac{[\alpha]}{\beta'}\text{induction}}{\exists x\beta'}\exists^{+}}{\exists x\beta'}\exists^{-}$$

```
... | yes refl with subIdentFunc α[x/y]≡β
...             | refl =
  close
   (assembled-context Γ⊢∃xα)
   (λ x₁ z z₁ → z z₁ (λ z₂ → z₂ (λ z₃ → z₃)))
   (existelim (all⟨ V↓ x β' ⟩ all∪ (all- all⟨ y∉α ⟩))
    Γ⊢∃xα
    (existintro (varterm x) x (ident β' x)
     (⟨→⟩ (renameIff β≈β')
      (assume α))))
```

Now, consider the other direction.

```
⟨←⟩ (renameIff {Γ} {V x α} {V y β'} (V/ y∉α α[x/y]≡β β≈β')) Γ⊢∃yβ' =
```

$$\dfrac{\begin{matrix}\Gamma \\ \vdots \\ \exists y\beta'\end{matrix} \qquad \dfrac{\dfrac{[\beta']}{\beta}\text{induction}}{\exists x\alpha}\exists^{+}}{\exists x\alpha}\exists^{-}$$

```
  close
   (assembled-context Γ⊢∃yβ')
   (λ x₁ x₂ x₃ → x₂ x₃ λ x₄ → x₄ λ x₅ → x₅)
   (existelim (all⟨ V x y∉α ⟩ all∪ (all- all⟨- [ refl ] ⟩))
    Γ⊢∃yβ'
    (existintro (varterm y) x α[x/y]≡β
     (⟨←⟩ (renameIff β≈β')
      (assume β'))))
```

The third case is the dual of the second.

```
⟨→⟩ (renameIff {Γ} {V x α} {V y β} (V/' α≈α' y∉α' α'[x/y]≡β')) Γ⊢∃xα
    with varEq x y
```

If $x = y$, then $\alpha' = \beta'$.

$$\dfrac{\begin{array}{c}\Gamma\\ \vdots\\ \exists x\alpha\end{array} \qquad \dfrac{\dfrac{[\alpha]}{\beta'}\,\text{induction}}{\exists x\beta'}\,\exists^+}{\exists x\beta'}\,\exists^-$$

```
... | yes refl rewrite subIdentFunc α'[x/y]≡β' =
  close
   (assembled-context Γ⊢∃xα)
   (λ x₁ z z₁ → z z₁ (λ z₂ → z₂ (λ z₃ → z₃)))
   (existelim (all⟨ V↓ x β ⟩ all∪ (all- all⟨- [ refl ] ⟩))
    Γ⊢∃xα
    (existintro (varterm x) x (ident β x)
     (⟨→⟩ (renameIff α≈α')
      (assume α))))
```

Otherwise, because $\alpha'[x/y] = \beta'$, we have $\alpha' = \beta'[y/x]$, and $x$ is not free in $\beta'$, and so is not free in $\exists y\beta'$.

$$\dfrac{\begin{array}{c}\Gamma\\ \vdots\\ \exists x\alpha\end{array} \qquad \dfrac{\dfrac{[\alpha]}{\alpha'}\,\text{induction}}{\exists y\beta'}\,\exists^+}{\exists y\beta'}\,\exists^-$$

```
... | no x≢y   =
  close
   (assembled-context Γ⊢∃xα)
   (λ x₁ z z₁ → z z₁ (λ z₂ → z₂ (λ z₃ → z₃)))
   (existelim (all⟨ V y (subNotFree (varterm x≢y) α'[x/y]≡β') ⟩
               all∪ (all- all⟨- [ refl ] ⟩))
    Γ⊢∃xα
    (existintro (varterm x) y (subInverse y∉α' α'[x/y]≡β')
     (⟨→⟩ (renameIff α≈α')
      (assume α))))
```

Consider the other direction.

```
⟨←⟩ (renameIff {Γ} {V x α} {V y β'} (V/' α≈α' y∉α' α'[x/y]≡β')) Γ⊢∃yβ' =
```

$$\dfrac{\begin{array}{c}\Gamma\\ \vdots\\ \exists y\beta'\end{array} \qquad \dfrac{\dfrac{[\beta']}{\exists x\alpha'}\,\exists^+ \qquad \dfrac{\dfrac{[\alpha']}{\alpha}\,\text{induction}}{\exists x\alpha}\,\exists^+}{\exists x\alpha}\,\exists^-}{\exists x\alpha}\,\exists^-$$

```
  close
   (assembled-context Γ⊢∃yβ')
   (λ x z z₁ → z z₁ (λ z₂ → z₂ (λ z₃ z₄ → z₄ z₃ (λ z₅ → z₅ (λ z₆ → z₆)))))
   (existelim (all⟨ V x (≈notFree (≈sym α≈α') y∉α') ⟩
               all∪ (all- (all⟨- [ refl ] ⟩ all∪ (all- all⟨ y∉α' ⟩))))
    Γ⊢∃yβ'
    (existelim (all⟨ V↓ x α ⟩ all∪ (all- all⟨- [ refl ] ⟩))
     (existintro (varterm y) x α'[x/y]≡β'
      (assume β'))
     (existintro (varterm x) x (ident α x)
      (⟨←⟩ (renameIff α≈α')
       (assume _)))))
```

□

We can conclude that examining formulae only on an intensional level does not restrict the deductive power of the system.

There is a dual structure in the proofs above, in the quantifier cases where the bound variable is renamed. Some proofs are straightforward in that they eliminate the quantifier, insert the derivation of the equivalent subcomponent by induction, then reintroduce the quantifier. Others are more complex, in that require an extra introduction and elimination step. The straightforward proofs are for the forward direction for `Λ/` and `V/'`, and the reverse direction for `Λ/'` and `V/`, while the complex proofs are the forward direction for `Λ/'` and `V/`, and the reverse direction for `Λ/` and `V/'`. Since the forward direction of each of these rules is the same as the reverse direction of its dual, we see that it would be simplest to do renaming with the rules `Λ/` and `V/'`, and have `Λ/'` and `V/` be the derived rules, if possible.

# 10 Scheme.lagda

The previous modules define the language of natural deduction. This system can be used to show that certain first-order formulae are derivable in minimal logic. It is common in logical enquiries to examine proofs regarding axiom schemes, as we will do later (see also [6, 12, 13]).

We define some metalanguage concepts. A *scheme* is often thought of as a formula containing schematic variables, which can be replaced by subformulae to produce a new formula. The following notion is more general than this; instead, a scheme is just constructed from a function from (a vector of) formulae to a formula.

```
record Scheme : Set where
  constructor scheme
  field
    arity : ℕ
    name  : String
    inst  : Vec Formula arity → Formula
```

Defining this as a type using a vector, instead of simply using functions, means that all schemes of all arities are collected under the same type (`Scheme`), which makes it possible to define a single function for typesetting scheme proofs later. The definition makes no restriction on the structure of the instances of the scheme, and is not able to put requirements on variable freedom.

A scheme is derivable if every instance of the scheme is derivable. A list $\Omega s$ of schemes is stronger than a scheme $\Phi$ if every instance of $\Phi$ is derivable from finitely many instances of schemes in $\Omega s$. Equivalently, $\Omega s$ is stronger than $\Phi$ if the derivability of $\Omega s$ implies the derivability of $\Phi$.

```
Derivable : Scheme → Set₁
Derivable S = ∀ αs → ⊢ (Scheme.inst S αs)

infix 1 _⊃_
_⊃_ : List Scheme → Scheme → Set₁
Ωs ⊃ Φ = (∀ ω → ω List.∈ Ωs → Derivable ω) → Derivable Φ
```

Because it is nicer to work with $n$-ary functions than unary functions taking $n$-ary vectors, we define the following notation for creating schemes from functions,

```
nullaryscheme : String → Formula → Scheme
unaryscheme   : String → (Formula → Formula) → Scheme
binaryscheme  : String → (Formula → Formula → Formula) → Scheme
```

```
nullaryscheme s f = scheme 0 s λ { [] → f }
unaryscheme   s f = scheme 1 s λ { (α :: []) → f α }
binaryscheme  s f = scheme 2 s λ { (α :: β :: []) → f α β }
```

expressing derivability for functions,

```
infix 1 ⊢₀_ ⊢₁_ ⊢₂_

⊢₀_ : Formula → Set₁
⊢₁_ : (Formula → Formula) → Set₁
⊢₂_ : (Formula → Formula → Formula) → Set₁

⊢₀ s =           ⊢ s
⊢₁ s = ∀ α   → ⊢ s α
⊢₂ s = ∀ α β → ⊢ s α β
```

and turning derivability of schemes into derivability of functions.

```
descheme₀ : {f : Vec Formula 0 → Formula}
          → (∀ as → ⊢ f as) → ⊢ f []
descheme₁ : {f : Vec Formula 1 → Formula}
          → (∀ as → ⊢ f as) → ∀ α → ⊢ f (α :: [])
descheme₂ : {f : Vec Formula 2 → Formula}
          → (∀ as → ⊢ f as) → ∀ α β → ⊢ f (α :: β :: [])

descheme₀ ⊢S     = ⊢S []
descheme₁ ⊢S α   = ⊢S (α :: [])
descheme₂ ⊢S α β = ⊢S (α :: (β :: []))
```

## 11 Example: the drinker paradox

We give an example of proving scheme derivability. We will also use a module for outputting natural deduction trees as LaTeX.

```
open import Texify
```

The code for this is entirely computational, and can be found in the appendix.

First, some syntactic sugar. The `pattern` notation causes Adga to recognise the notation in places where their values would be used in pattern matching, and moreover will use the notation in proofs created by proof search. Note that we are no longer using ⊥ and ¬ as defined previously for decidable predicates in the metalanguage; here they are in the language of formulae.

```
pattern ⊥ = atom (rel zero zero) []

pattern ¬  α = α ⇒ ⊥
pattern ¬¬ α = ¬ (¬ α)
```

Fix some variables.

```
pattern xvar = var zero
pattern yvar = var (suc zero)

x y : Term
x = varterm xvar
```

```
y = varterm yvar

pattern  ∀x  Φ = Λ xvar Φ
pattern  ∃x  Φ = V xvar Φ
pattern ¬∀x  Φ = ¬(∀x Φ)
pattern ¬∃x  Φ = ¬(∃x Φ)
pattern  ∀x¬ Φ = ∀x (¬ Φ)
pattern  ∃x¬ Φ = ∃x (¬ Φ)
```

Define a nullary and a unary predicate (in the language of formulae), which will be used to instantiate the scheme proofs for output as proof trees in LaTeX.

```
pattern Arel = rel 1 0
pattern A    = atom Arel []

pattern Prel = rel 5 1
pattern P t  = atom Prel (t :: [])
```

The indices used for $x$, $y$, $\perp$, $A$, and $P$ are arbitrary, but correspond to those used internally by the texify module, so they will be outputted with the appropriate names.

Define the schemes DNE (double negation elimination), EFQ (ex falso quodlibet), DP (the drinker paradox), and H$\epsilon$ (the dual of the drinker paradox). The latter two schemes will be described and examined in more detail in the next chapter.

```
dne efq dp hε : Formula → Formula
dne Φ   = ¬¬ Φ ⇒ Φ
efq Φ   = ⊥ ⇒ Φ
dp  Φx  = ∃x(Φx ⇒ ∀x Φx)
hε  Φx  = ∃x(∃x Φx ⇒ Φx)

DNE EFQ DP Hε : Scheme
DNE = unaryscheme "DNE"            dne
EFQ = unaryscheme "EFQ"            efq
DP  = unaryscheme "DP"             dp
Hε  = unaryscheme "H$\\epsilon$" hε
```

The natural deduction system used to define _⊢_ is for minimal logic. This can be extended to classical logic with the classical $\perp$ rule.

```
⊥c-rule : Set₁
⊥c-rule = ∀{Γ} → ∀ α
          →      Γ ⊢ ⊥
            --------------- ⊥c
          →  Γ - (¬ α) ⊢ α
```

Similarly, the intuitionistic $\perp$ rule

```
⊥i-rule : Set₁
⊥i-rule = ∀{Γ} → ∀ α
          →      Γ ⊢ ⊥
                ------- ⊥i
          →      Γ ⊢ α
```

gives an extension to intuitionistic logic.

**Proposition 11.0.1.** *The classical bottom rule holds if and only if DNE is derivable.*

*Proof.*

```
dne→⊥c-rule : ⊢₁ dne → ⊥c-rule
```

$$\dfrac{\dfrac{}{\neg\neg\alpha \to \alpha}\,\text{DNE} \qquad \dfrac{\begin{matrix}\Gamma,\ [\neg\alpha]\\ \vdots \\ \bot\end{matrix}}{\neg\neg\alpha}\,\to^{+}}{\alpha}\,\to^{-}$$

```
dne→⊥c-rule ⊢dne α Γ⊢⊥ = close
                          (assembled-context (arrowintro (¬ α) Γ⊢⊥))
                          (λ x₁ z₁ z₂
                           → z₂ (λ z₃ → z₁ (λ z₄ → z₄) (λ z₄ → z₄ z₃)))
                          (arrowelim
                           (⊢dne α)
                           (arrowintro (¬ α)
                            Γ⊢⊥))
```

```
⊥c-rule→dne : ⊥c-rule → ⊢₁ dne
```

$$\dfrac{\dfrac{\dfrac{[\neg\neg\alpha] \qquad [\neg\alpha]}{\bot}\,\to^{-}}{\alpha}\,\bot_c}{\neg\neg\alpha \to \alpha}\,\to^{+}$$

```
⊥c-rule→dne ⊢⊥c-rule α = close
                          from∅
                          (λ x₁ z₁ z₂
                           → z₂ (z₁ (λ z₃ z₄ → z₄ (λ z₅ z₆ → z₆ z₃ z₅))))
                          (arrowintro (¬¬ α)
                           (⊢⊥c-rule α
                            (arrowelim
                             (assume (¬¬ α))
                             (assume (¬ α)))))
```

□

**Proposition 11.0.2.** *The intuitionistic bottom rule holds if and only if EFQ is derivable.*

*Proof.*

```
efq→⊥i-rule : ⊢₁ efq → ⊥i-rule
```

$$\dfrac{\dfrac{}{\bot \to \alpha}\,\text{EFQ} \qquad \begin{matrix}\Gamma\\ \vdots \\ \bot\end{matrix}}{\alpha}\,\to^{-}$$

```
efq→⊥i-rule ⊢efq α Γ⊢⊥ = close
                          (assembled-context Γ⊢⊥)
                          (λ x₁ z₁ → z₁ (λ z₂ → z₂))
                          (arrowelim
                           (⊢efq α)
                          Γ⊢⊥)
```

```
⊥i-rule→dne : ⊥i-rule → ⊢₁ efq
```

$$\dfrac{\dfrac{[\bot]}{\alpha}\ \bot_i}{\bot \to \alpha}\ \to^+$$

```
⊥i-rule→dne ⊢⊥i-rule α = close
                          from∅
                          (λ x1 z1 z2 → z2 (z1 (λ z3 → z3)))
                          (arrowintro ⊥
                           (⊢⊥i-rule α
                            (assume ⊥)))
```

∎

**Proposition 11.0.3.** *DP holds in classical logic.*

*Proof.* We show that if DNE is derivable then DP is derivable, meaning that DP is weaker than DNE. For illustrative purposes, lines given by Agda's proof search are marked with `{- Auto -}` in the next proof. The remainder of the proof, with the exception of the `close` function call, corresponds exactly to doing natural deduction by hand, from the bottom up. As the proof tree is developed, Agda displays the subgoal of each hole in the deduction, and will accept valid subproofs and formulae only. In this way, Agda not only verifies the deduction after it has been completed, but also acts as a proof assistant for natural deduction.

```
dne→dp : ⊢1 dne → ⊢1 dp
dne→dp ⊢dne α = close
  {- Auto -}      from∅
  {- Auto -}      (λ x1 z1 z2 → z2 (z1 (λ z3 → z3) (λ z3 → z3 (λ z4 z5 → z5 z4
  {- Auto -}       (λ z6 → z6 (λ _ z7 → z7 (λ z8 → z8) (λ z8 → z8 (λ z9 z10
  {- Auto -}        → z10 z4 (λ z11 → z11 (λ z12 z13 → z13 (λ z14 → z14)
  {- Auto -}         (λ z14 → z14 (λ _ z15 → z15 z9 z12))))))))))))
                  (arrowelim
                   (⊢dne (dp α))
                   (arrowintro (¬ (dp α))
                    (arrowelim
                     (assume (¬ (dp α)))
                     (existintro x xvar
  {- Auto -}          (ident (α ⇒ ∀x α) xvar)
                      (arrowintro α
                       (univintro xvar
  {- Auto -}             (all∅ all∪ (all- (all⟨ V↓ xvar (α ⇒ ∀x α) ⇒ atom [] ⟩
  {- Auto -}              all∪ (all- (all∅ all∪ (all- (all⟨- ¬∀x α :: (α ::
  {- Auto -}               [ refl ]) ⟩ all∪ all⟨- ¬∀x α :: [ refl ] ⟩)))))))
                        (arrowelim
                         (⊢dne α)
                         (arrowintro (¬ α)
                          (arrowelim
                           (assume (¬ (dp α)))
                           (existintro x xvar
  {- Auto -}                (ident (α ⇒ ∀x α) xvar)
                            (arrowintro α
                             (arrowelim
                              (⊢dne (∀x α))
                              (arrowintro (¬∀x α)
                               (arrowelim
                                (assume (¬ α))
                                (assume α)))))))))))))))
```

∎

The above is a general derivation of an arbitrary instance of DP using instances of DNE. We use this proof to construct the scheme relation '⊃', for outputting as LaTeX.

```
DNE⊃DP : DNE :: [] ⊃ DP
DNE⊃DP ⊢lhs (α :: []) = dne→dp (descheme₁ (⊢lhs DNE [ refl ])) α
dp-prooftree = texreduce DP (P x :: []) DNE⊃DP
```

The final line gets the deduction tree for the instance DP($Px$), which is shown below, with instances of DP abbreviated, and split into two, due to page constraints.

$$\dfrac{[\neg\exists x\,(Px \to \forall x Px)] \qquad \dfrac{\dfrac{\dfrac{\dfrac{}{\neg\neg\forall x Px \to \forall x Px}\,\text{DNE} \qquad \dfrac{\dfrac{[\neg Px] \qquad [Px]}{\bot}\,{\to^-}}{\neg\neg\forall x Px}\,{\to^+}}{\forall x Px}\,{\to^-}}{Px \to \forall x Px}\,{\to^+}}{\exists x\,(Px \to \forall x Px)}\,{\exists^+}}{\dfrac{\bot}{\vdots}}\,{\to^-}$$

$$\dfrac{\dfrac{}{\neg\neg\mathrm{DP}(Px) \to \mathrm{DP}(Px)}\,\text{DNE} \qquad \dfrac{\dfrac{[\neg\exists x\,(Px \to \forall x Px)] \qquad \dfrac{\dfrac{\dfrac{\dfrac{\dfrac{}{\neg\neg Px \to Px}\,\text{DNE} \qquad \dfrac{\dfrac{\vdots}{\bot}}{\neg\neg Px}\,{\to^+}}{Px}\,{\to^-}}{\forall x Px}\,{\forall^+}}{Px \to \forall x Px}\,{\to^+}}{\exists x\,(Px \to \forall x Px)}\,{\exists^+}}{\bot}\,{\to^-}}{\neg\neg\exists x\,(Px \to \forall x Px)}\,{\to^+}}{\exists x\,(Px \to \forall x Px)}\,{\to^-}$$

**Proposition 11.0.4.** *The dual of the drinker paradox also holds in classical logic.*

*Proof.*

```
dne→hε : ⊢₁ dne → ⊢₁ hε
dne→hε ⊢dne α = close
                 from∅
                 (λ x₁ z₁ z₂ → z₂ (z₁ (λ z₃ → z₃) (λ z₃ → z₃ (λ z₄ z₅ → z₅ z₄
                  (λ z₆ → z₆ (λ z₇ z₈ → z₈ (λ z₉ → z₉) (λ z₉ → z₉ (λ _ z₁₀
                   → z₁₀ z₇ (λ z₁₁ → z₁₁ (λ z₁₂ z₁₃ → z₁₃ z₄ (λ z₁₄ → z₁₄
                    (λ _ → z₁₂)))))))))))))
                 (arrowelim
                  (⊢dne (hε α))
                  (arrowintro (¬ (hε α))
                   (arrowelim
                    (assume (¬ (hε α)))
                    (existintro x xvar (ident (∃x α ⇒ α) xvar)
                     (arrowintro (∃x α)
                      (arrowelim
                       (⊢dne α)
                       (arrowintro (¬ α)
                        (existelim
```

```
                          (all⟨ atom [] ⟩ all∪ (all- (all⟨ V↓ xvar (∃x α ⇒ α)
                           ⇒ atom [] ⟩ all∪ (all- all⟨- ∃x α :: [ refl ] ⟩))))
                          (assume (∃x α))
                          (arrowelim
                           (assume (¬ (hε α)))
                           (existintro x xvar (ident (∃x α ⇒ α) xvar)
                            (arrowintro (∃x α)
                             (assume α))))))))))))
```
□

We extract the proof tree for $\mathrm{H}\epsilon(Px)$.

```
DNE⊃Hε : DNE :: [] ⊃ Hε
DNE⊃Hε ⊢lhs (α :: []) = dne→hε (descheme₁ (⊢lhs DNE [ refl ])) α
hε-prooftree = texreduce Hε (P x :: []) DNE⊃Hε
```

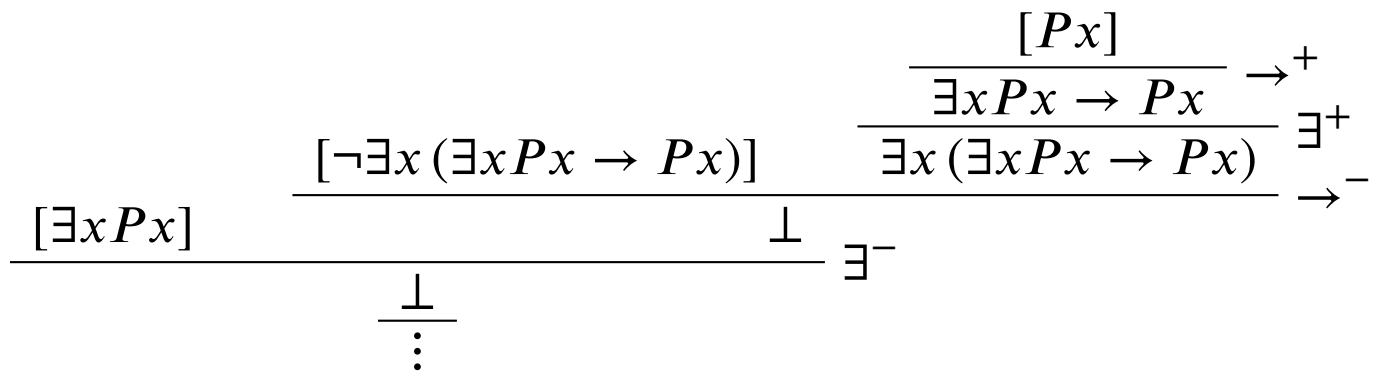

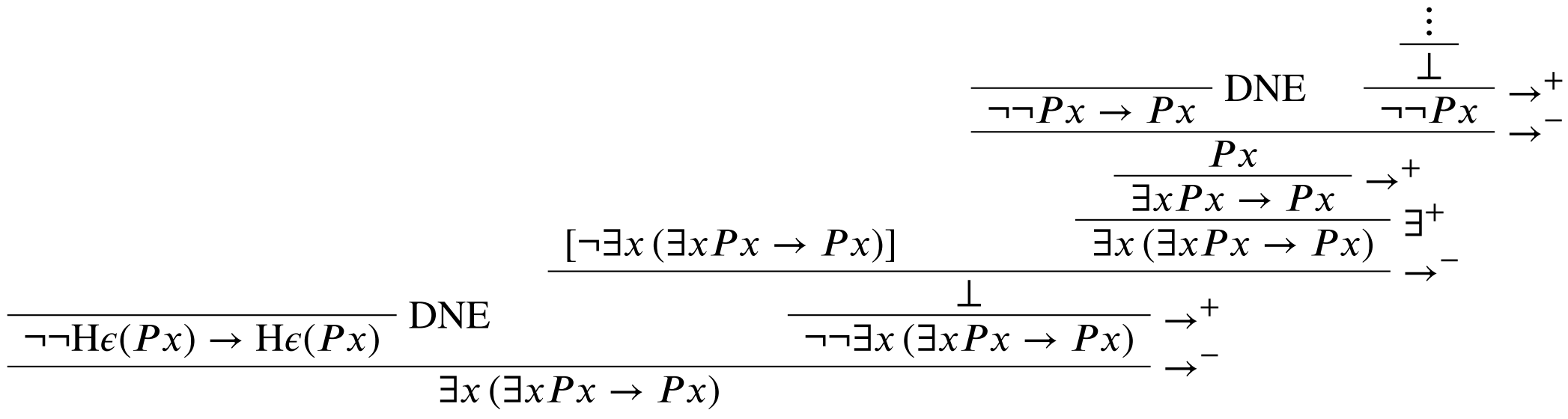

As a final example, consider the law of excluded middle, and a general form of the limited principle of omniscience.[9]

```
lem glpo : Formula → Formula
lem  Φ = Φ ∨ (¬ Φ)
glpo Φ = ∀x (¬ Φ) ∨ ∃x Φ

LEM GLPO : Scheme
LEM  = unaryscheme "LEM"  lem
GLPO = unaryscheme "GLPO" glpo
```

Recall that equivalent formulae are equivalently derivable, so from GLPO we may derive a form with any other quantifying variable. Therefore while the variable $x$ is fixed, it can be expected that LEM and GPO are equivalent with respect to derivability. That is, in an extension of minimal logic where one is derivable, the other should also be derivable. The former leads to the latter in a straightforward manner. The other direction is more complicated, since $\Phi$ could have $x$ free.

We show first that when deriving LEM($\Phi$), we may assume without loss of generality that $x$ is not free in $\Phi$, by showing that if LEM is derivable in this restricted case then it is derivable in general.

[9] This is general in the sense that it is not over a binary sequence, like that of [3], but rather over a predicate which may not be decidable.

*Proof.* Given any formula $\alpha$, there is a fresh variable $\omega$ which appears nowhere in $\alpha$ and which differs from $x$. Then $\alpha[x/\omega]$ exists, with $x$ not free, and $\alpha[x/\omega][\omega/x] = \alpha$. Now if LEM holds for $\alpha[x/\omega]$ then it holds for $\alpha$, by the following proof tree.

$$\dfrac{\dfrac{\overline{\alpha[x/\omega] \vee \neg\alpha[x/\omega]}}{\forall\omega\,(\alpha[x/\omega] \vee \neg\alpha[x/\omega])}\,\forall^+}{\alpha \vee \neg\alpha}\,\forall^-$$

Hence we may derive LEM by deriving it only for formulae for which $x$ is not free. This is formalised in Agda as follows.

```
wlog-lem : (∀ α → xvar NotFreeIn α → ⊢ (lem α)) → ⊢₁ lem
wlog-lem ⊢nflem α = close
                     from∅
                     (λ x₁ z₁ z₂ → z₂ z₁)
                     (univelim x lemαω[ω/x]≡lemα
                      (univintro ωvar all∅
                       (⊢nflem αω x∉αω)))
  where
```

Compute the fresh variable, and use its construction to get that it is fresh in $\alpha$ and not equal to $x$.

```
    ω,ωFresh,x≢ω : Σ Variable (λ ω → Σ (ω FreshIn α) (λ _ → xvar ≢ ω))
    ω,ωFresh,x≢ω with fresh (∀x α)
    ...            | ω , Λ x≢ω ωFrα = ω , ωFrα , x≢ω
```

We therefore have a variable $\omega$ which is not free in $\alpha$, which is free for $x$ in $\alpha$, and which differs from $x$.

```
    ωvar           : Variable
    ω∉α            : ωvar NotFreeIn α
    ωFreeForxInα   : (varterm ωvar) FreeFor xvar In α
    x≢ω            : xvar ≢ ωvar
    ωvar           = fst ω,ωFresh,x≢ω
    ω∉α            = freshNotFree (fst (snd ω,ωFresh,x≢ω))
    ωFreeForxInα   = freshFreeFor (fst (snd ω,ωFresh,x≢ω)) xvar
    x≢ω            = snd (snd ω,ωFresh,x≢ω)
```

Now, compute $\alpha_\omega = \alpha[x/\omega]$.

```
    αω         : Formula
    α[x/ω]≡αω : α [ xvar / _ ]≡ αω
    αω         = fst (α [ xvar / ωFreeForxInα ])
    α[x/ω]≡αω = snd (α [ xvar / ωFreeForxInα ])
```

By the construction of $\omega$, the substitution is reversible, so $\mathrm{LEM}(\alpha_\omega)[\omega/x] = \mathrm{LEM}(\alpha)$.

```
    lemαω[ω/x]≡lemα : (lem αω) [ ωvar / _ ]≡ (lem α)
    lemαω[ω/x]≡lemα = subInverse
                         (ω∉α ∨ (ω∉α ⇒ atom []))
                         (α[x/ω]≡αω ∨ (α[x/ω]≡αω ⇒ notfree (atom [])))
```

Finally, $x$ will not be free after it has been substituted out of $\alpha$.

```
    x∉αω : xvar NotFreeIn αω
    x∉αω = subNotFree (varterm x≢ω) α[x/ω]≡αω
```

□

We can now show that GLPO is stronger than LEM, without worrying about the quantifier variable.

```
glpo→xnf→lem : ⊢₁ glpo → ∀ α → xvar NotFreeIn α → ⊢ (lem α)
glpo→xnf→lem ⊢glpo α x∉α = close
                              from∅
                              (λ x₁ z₁ z₂ → z₂ (z₁ (λ z₃ → z₃) (λ z₃ → z₃)
                               (λ z₄ → z₄ (λ z₅ → z₅)) (λ z₄ → z₄ (λ z₅ z₆ → z₆
                                z₅ (λ z₇ → z₇ (λ z₈ → z₈)))))))
                              (disjelim
                               (⊢glpo α)
                               (disjintro₂ α
                                (univelim x (ident (¬ α) xvar)
                                 (assume (∀x¬ α))))
                               (disjintro₁ (¬ α)
                                (existelim (all⟨ x∉α ⟩ all∪ (all- all⟨ x∉α ⟩))
                                 (assume (∃x α))
                                 (assume α))))
```

Now, LEM can be obtained directly from GLPO. The proof tree for the restricted form of LEM is inserted into the proof tree from `wlog-lem`.

```
glpo→lem : ⊢₁ glpo → ⊢₁ lem
glpo→lem ⊢glpo = wlog-lem (glpo→xnf→lem ⊢glpo)

GLPO⊃LEM : GLPO :: [] ⊃ LEM
GLPO⊃LEM ⊢lhs (α :: []) = glpo→lem (descheme₁ (⊢lhs GLPO [ refl ])) α
```

No computation of a fresh variable has occurred yet, since the variable depends on the instance of LEM we want to derive. Extracting the proof tree for LEM($Px$), the `fresh` function computes that $y$ is fresh, and so the proof tree below is produced.

```
glpo→lem-prooftree = texreduce LEM (P x :: []) GLPO⊃LEM
```

$$
\dfrac{\dfrac{\dfrac{\dfrac{}{\forall x \neg Py \lor \exists x Py}\,\text{GLPO} \qquad \dfrac{\dfrac{[\forall x \neg Py]}{\neg Py}\,\forall^-}{Py \lor \neg Py}\,\lor^+ \qquad \dfrac{\dfrac{[\exists x Py] \qquad [Py]}{Py}\,\exists^-}{Py \lor \neg Py}\,\lor^+}{Py \lor \neg Py}\,\lor^-}{\forall y\,(Py \lor \neg Py)}\,\forall^+}{Px \lor \neg Px}\,\forall^-
$$

# Appendices

## A Context as a list

### A.1 Computational definition

Lists could be used for the context of natural deduction trees, instead of using ensembles. The operations for removal and union are still needed.

```
infixl 5 _-_
_-_ : List Formula → Formula → List Formula
[] - β = []
(α :: αs) - β with formulaEq α β
((β :: αs) - .β) | yes refl = αs - β
((α :: αs) -  β) | no  _    = α :: (αs - β)

infixr 5 _∪_
_∪_ : List Formula → List Formula → List Formula
[]        ∪ ys = ys
(x :: xs) ∪ ys = x :: (xs ∪ ys)
```

However, it is now more complicated to prove that a given deduction's context is a subset of the permitted open assumptions. It is necessary to reason about the result of a computation. Begin with the following trivial lemma.

```
eqcontext : ∀{α Δ Γ}→ Δ ≡ Γ → Δ ⊢ α → Γ ⊢ α
eqcontext refl x = x
```

The proof for $\vdash \alpha \to \alpha$ is as follows

```
arrow-example : ∀ α → ⊢ α ⇒ α
arrow-example α = eqcontext closed
                    (arrowintro α
                     (assume α))
                   where
                     closed : ((α :: []) - α) ≡ []
                     closed with formulaEq α α
                     ...      | yes refl = refl
                     ...      | no  α≢α  = ⊥-elim (α≢α refl)
```

To examine what the result of `(α :: []) - α` is, we must examine the pattern matching that occurs on the result of `formulaEq α α`. In the real case where `α ≡ α` holds, $\alpha$ is removed from the list, and the proof is `refl`. However, we must also consider the case where `α ≢ α` (which we prove using absurdity).

Now, consider a more complicated example; we prove that $\alpha \to \beta \to \gamma \vdash \beta \to \alpha \to \gamma$.

```
reorder : ∀ α β γ → α ⇒ β ⇒ γ :: [] ⊢ β ⇒ α ⇒ γ
reorder α β γ = eqcontext closed
                 (arrowintro β
                  (arrowintro α
                   (arrowelim
                    (arrowelim
                     (assume (α ⇒ β ⇒ γ))
                     (assume α))
                    (assume β))))
  where
    closed : ((α ⇒ β ⇒ γ :: α :: β :: []) - α - β) ≡ α ⇒ β ⇒ γ :: []
    closed with formulaEq (α ⇒ β ⇒ γ) α
    closed | yes ()
    closed | no _ with formulaEq α α
    closed | no _ | no  α≢α  = ⊥-elim (α≢α refl)
    closed | no _ | yes refl with formulaEq β α
    closed | no _ | yes refl | yes refl with formulaEq (α ⇒ β ⇒ γ) β
    closed | no _ | yes refl | yes refl | yes ()
    closed | no _ | yes refl | yes refl | no _ = refl
    closed | no _ | yes refl | no _ with formulaEq (α ⇒ β ⇒ γ) β
    closed | no _ | yes refl | no _ | yes ()
    closed | no _ | yes refl | no _ | no _ with formulaEq β β
    closed | no _ | yes refl | no _ | no _ | yes refl = refl
    closed | no _ | yes refl | no _ | no _ | no  β≢β  = ⊥-elim (β≢β refl)
```

Each equality check must be examined. Clearly this becomes unwieldy, even in simple cases. Moreover, Agda's proof search will not create `with` blocks, and so is of little use here.

## A.2 Expanded context definition

In a similar fashion to [5], we could define a proof system which is similar to natural deduction, which does not use list computation in the main deduction rules. Instead, include a deduction rule for weakening the context on the left, and allow `assume` to weaken the context on the right. We show only the rules for implication and universal generalisation.

```
_++_ : List Formula → List Formula → List Formula
[]        ++ ys = ys
(x :: xs) ++ ys = x :: (xs ++ ys)

infix 1 _⊢_ ⊢_
data _⊢_ : List Formula → Formula → Set where

  assume      : ∀{Γ} → (α : Formula)
              →                                α :: Γ ⊢ α

  weaken      : ∀{Γ α} → (Δ : List Formula)
              →                                  Γ ⊢ α
                                              --------------
              →                                (Δ ++ Γ) ⊢ α

  arrowintro : ∀{Γ β} → (α : Formula)
              →                                α :: Γ ⊢ β
                                               ----------- ⇒⁺
              →                                Γ ⊢ α ⇒ β

  arrowelim   : ∀{Γ α β}
              →                        Γ ⊢ α ⇒ β     →     Γ ⊢ α
                                      ------------------------- ⇒⁻
              →                                 Γ ⊢ β

  univintro   : ∀{Γ α} → (x : Variable)
              → x NotFreeInAll Γ
              →                                  Γ ⊢ α
                                               ----------- ∀⁺
              →                                Γ ⊢ Λ x α

  univelim    : ∀{Γ α x α[x/r]} → (r : Term)
              → α [ x / r ]≡ α[x/r]
              →                                Γ ⊢ Λ x α
                                               ------------ ∀⁻
              →                                Γ ⊢ α[x/r]
```

This system does not describe natural deduction, since the context is not the same as it is for natural deduction. Extra formulae are assumed. It also requires weakening at each assumption. Weakening could be built into the assumption definition, and can be solved by proof search, but it is not a usual consideration when doing natural deduction by hand.

This system works for propositional logic. We again prove that $\alpha \to \beta \to \gamma \vdash \beta \to \alpha \to \gamma$.

```
reorder : ∀ α β γ → α ⇒ β ⇒ γ :: [] ⊢ β ⇒ α ⇒ γ
reorder α β γ = arrowintro β
                 (arrowintro α
                  (arrowelim
                   (arrowelim
                    (weaken (α :: β :: [])
                     (assume (α ⇒ β ⇒ γ)))
                     (assume α))
```

```
                    (weaken (α :: [])
                     (assume β))))
```

The added assumptions become an issue for the first order case, due to the restrictions on free variables. Consider the following proof.

$$\dfrac{\forall x Qx \to Px \qquad \dfrac{\dfrac{\dfrac{\forall x\,(\forall x A \to Qx)}{\forall x A \to Qx}\,\forall^- \qquad \forall x A}{Qx}\,\to^-}{\forall x Qx}\,\forall^+}{Px}\,\to^-$$

This is a valid natural deduction, and it was checked with the ensemble-based natural deduction system. However, this proof tree does not satisfy the above rules, since $\forall x A \to Qx$ would have to be made an extra assumption above the deduction of $Qx$ by weakening. This means that the universal generalisation introduction is not valid, since $x$ is free in an open assumption.

## B Texify.lagda

The following code defines a function called `texdeduction` for outputting proof trees as LaTeX, using the *bussproofs* package. The function `texreduce` does the same for derivations from schemes to schemes.

```
module Texify where

open import Agda.Builtin.Bool
open import Agda.Builtin.Nat renaming (Nat to ℕ)
open import Agda.Builtin.String

open import Decidable hiding (⊥ ; ¬_)
open import Deduction
open import Ensemble
open import Formula
open import List
open import Scheme
open import Vec
open import Sugar

TEXarrowintro = "$\\rightarrow^+$"
TEXarrowelim  = "$\\rightarrow^-$"
TEXconjintro  = "$\\land^+$"
TEXconjelim   = "$\\land^-$"
TEXdisjintro  = "$\\lor^+$"
TEXdisjelim   = "$\\lor^-$"
TEXunivintro  = "$\\forall^+$"
TEXunivelim   = "$\\forall^-$"
TEXexistintro = "$\\exists^+$"
TEXexistelim  = "$\\exists^-$"
TEXarrow      = " \\rightarrow "
TEXand        = " \\land "
TEXor         = " \\lor "
TEXforall     = "\\forall"
TEXexists     = "\\exists"
TEXnot        = "\\lnot"
```

```
TEXbot         = "\\bot"

lp = "\\left("
rp = "\\right)"


-- String manipulation
_>>_ = primStringAppend
infixr 1 _>>_

wrap : String → String
wrap s = "{" >> s >> "}"

-- Instead of using stdlib
strnum : ℕ → String
strnum zero = "0"
strnum (suc n) = wrap ("s(" >> strnum n >> ")")

strrel : Relation → String
strrel (rel 0 k) = "\\bot"
strrel (rel 1 k) = "A"
strrel (rel 2 k) = "B"
strrel (rel 3 k) = "C"
strrel (rel 4 k) = "D"
strrel (rel 5 k) = "P"
strrel (rel 6 k) = "Q"
strrel (rel (suc (suc (suc (suc (suc (suc n)))))) k) = "R_" >> strnum n

strvar : Variable → String
strvar xvar = "x"
strvar yvar = "y"
strvar zvar = "z"
strvar (var n) = "v_" >> strnum n

-- The constants are the natural numbers
strfunc : Function → String
strfunc (func n k) = "f_" >> strnum n

join : String → List String → String
join delim [] = ""
join delim (s :: []) = s
join delim (s :: ss@(_ :: _)) = s >> delim >> join delim ss

joinmap : {A : Set} → String → (A → String) → List A → String
joinmap delim f [] = ""
joinmap delim f (x :: []) = f x
joinmap delim f (x :: xs@(_ :: _)) = f x >> delim >> joinmap delim f xs

texterm : Term → String
textermvec : ∀{n} → Vec Term n → String

texterm (varterm x) = wrap (strvar x)
texterm t0 = wrap "0"
```

```
texterm t1 = wrap "1"
texterm t2 = wrap "2"
texterm t3 = wrap "3"
texterm t4 = wrap "4"
texterm t5 = wrap "5"
texterm (functerm (func n f) ts) with n
...     | zero   = wrap (strfunc (func n f))
...     | suc _ = wrap (strfunc (func n f) >> lp >> textermvec ts >> rp)

textermvec [] = ""
textermvec (t :: []) = texterm t
textermvec (t :: ts@(_ :: _)) = texterm t >> ", " >> textermvec ts


texformula : Formula → String

parenformula : Formula → String
parenformula p@(atom _ _) = texformula p
parenformula p@(_ ⇒ b) with formulaEq b ⊥ | texformula p
...                     | yes _ | tfp = tfp
...                     | no _  | tfp = lp >> tfp >> rp
parenformula p@(_ ∧ _) = lp >> texformula p >> rp
parenformula p@(_ ∨ _) = lp >> texformula p >> rp
parenformula p@(Λ _ _) = texformula p
parenformula p@(V _ _) = texformula p

texformula a@(atom f ts) with formulaEq a ⊥
...                             | yes _ = TEXbot
texformula (atom (rel n k) ts) | no _  with k
...     | zero      = strrel (rel n k)
...     | suc zero = strrel (rel n k) >> textermvec ts
texformula (atom (rel n k) (x :: y :: [])) | no _
        | suc (suc zero) = texterm x >> strrel (rel n k) >> texterm y
...     | suc (suc _)     = strrel (rel n k)  >> lp >> textermvec ts >> rp
texformula (a ⇒ b) with formulaEq b ⊥
...     | yes _  = TEXnot >> wrap (parenformula a)
...     | no _  = parenformula a >> TEXarrow >> parenformula b
texformula (a ∧ b) = parenformula a >> TEXand >> parenformula b
texformula (a ∨ b) = parenformula a >> TEXor >> parenformula b
texformula (Λ x a) = TEXforall >> wrap(strvar x) >> parenformula a
texformula (V x a) = TEXexists >> wrap(strvar x) >> parenformula a

texformulae : List Formula → String
texformulae forms = joinmap ", " texformula forms


data Textree : Set where
  schemeax    : Formula → String → Textree
  openax      : Formula → Textree
  closedax    : Formula → Textree
  unaryinf    : Formula → String → Textree → Textree
  binaryinf   : Formula → String → Textree → Textree → Textree
  trinaryinf  : Formula → String → Textree → Textree → Textree → Textree
```

```
line : ℕ → String → String
line zero s = s >> "\n"
line (suc n) s = "\t" >> line n s

tag : String → String → String
tag f s = "\\" >> f >> "{" >> s >> "}"

label : ℕ → String → String
label i s = line i (tag "RightLabel" s)

inf : ℕ → String → Formula → String
inf i s x = line i (tag s ("$" >> (texformula x) >> "$"))

dis : ℕ → String → Formula → String
dis i s x = line i (tag s ("$\\left[" >> (texformula x) >> "\\right]$"))

texifytree : ℕ → Textree → String
texifytree i (schemeax x s)            = line i ("\\AxiomC{}")
                                         >> label i s
                                         >> inf i "UnaryInfC" x
texifytree i (openax x)                = inf i "AxiomC" x
texifytree i (closedax x)              = dis i "AxiomC" x
texifytree i (unaryinf x s T)          = texifytree i T
                                         >> label i s
                                         >> inf i "UnaryInfC" x
texifytree i (binaryinf x s T₁ T₂)      = texifytree i T₁
                                         >> texifytree (i + 1) T₂
                                         >> label i s
                                         >> inf i "BinaryInfC" x
texifytree i (trinaryinf x s T₁ T₂ T₃)  = texifytree i T₁
                                         >> texifytree (i + 1) T₂
                                         >> texifytree (i + 2) T₃
                                         >> label i s
                                         >> inf i "TrinaryInfC" x


dtot : ∀{α Γ} {ω : Ensemble Formula}
       → Assembled formulaEq ω → Γ ⊢ α → Textree
dtot {α} o (cite s d)            = schemeax α s
dtot {α} o (assume a) with Ensemble.dec∈ a o
...                     | yes _  = openax    α
...                     | no  _  = closedax  α
dtot {α} o (arrowintro a d)      = unaryinf  α TEXarrowintro (dtot o d)
dtot {α} o (arrowelim d₁ d₂)     = binaryinf α TEXarrowelim  (dtot o d₁)
                                                             (dtot o d₂)
dtot {α} o (conjintro d₁ d₂)     = binaryinf α TEXconjintro  (dtot o d₁)
                                                             (dtot o d₂)
dtot {α} o (conjelim d₁ d₂)      = binaryinf α TEXconjelim   (dtot o d₁)
                                                             (dtot o d₂)
dtot {α} o (disjintro₁ b d)      = unaryinf  α TEXdisjintro  (dtot o d)
dtot {α} o (disjintro₂ a d)      = unaryinf  α TEXdisjintro  (dtot o d)
```

```
dtot {α} o (disjelim d₁ d₂ d₃)   = trinaryinf α TEXdisjelim   (dtot o d₁)
                                                              (dtot o d₂)
                                                              (dtot o d₃)
dtot {α} o (univintro x _ d)     = unaryinf   α TEXunivintro  (dtot o d)
dtot {α} o (univelim r _ d)      = unaryinf   α TEXunivelim   (dtot o d)
dtot {α} o (existintro r x _ d)  = unaryinf   α TEXexistintro (dtot o d)
dtot {α} o (existelim _ d₁ d₂)   = binaryinf  α TEXexistelim  (dtot o d₁)
                                                              (dtot o d₂)
dtot {α} o (close _ _ d)         = dtot o d


texdeduction : ∀{Γ α} → Γ ⊢ α → String
texdeduction d = "\\begin{prooftree}\n"
                 >> texifytree 0 (dtot (assembled-context d) d)
                 >> "\\end{prooftree}\n"


-- We postulate that every instance of the stronger schemes is derivable. By
-- using cite, the deductions for these become irrelevant, so a string is
-- still produced. Postulates are not safe, but here this should not cause
-- problems, since this can only be used to produce strings.

-- The scheme y cannot be found implicitly because of how strings are defined.
texreduce : {xs : List Scheme}
            → (y : Scheme) → Vec Formula (Scheme.arity y) → xs ⊃ y → String
texreduce {xs} y αs xs⊃y = texdeduction (xs⊃y ⊢xs αs)
  where
    ⊢xs : (x : Scheme) → x List.∈ xs → Derivable x
    ⊢xs (scheme n name f) _ αs = cite name Ω
      where
        postulate Ω : ∅ ⊢ f αs

texprop : {xs : List Scheme}
            → (y : Scheme) → Vec Formula (Scheme.arity y) → xs ⊃ y → String
texprop {xs} y αs xs⊃y
  = "\\begin{proposition}\n"
    >> "$\\text{" >> joinmap "," Scheme.name xs >> "}"
       >> " \\supset \\text{" >> Scheme.name y >> "}$\n"
    >> "\\end{proposition}\n"
    >> "\\begin{proof}\n"
    >> "$ $\n"
    >> "\\vspace{-\\baselineskip}\n"
    >> (texreduce y αs xs⊃y)
    >> "\\vspace{-\\baselineskip}\n"
    >> "\\end{proof}\n"
```